\documentstyle[12pt]{article}
\begin{document}
\title{Quasi-invariant and pseudo-differentiable measures on
a non-Archimedean Banach space.II. Measures with values in
non-Archimedean fields.}
\author{Sergey V. Ludkovsky.}
\date{25.05.2001}
\maketitle
\begin{abstract}
Quasi-invariant and pseudo-differentiable measures on a
Banach space $X$ over a non-Archimedean locally compact 
infinite field with a non-trivial valuation
are defined and constructed. Measures are considered with values
in non-Archimedean fields, for example, the field $\bf Q_p$ of $p$-adic 
numbers. Theorems and criteria are formulated and
proved about
quasi-invariance and pseudo-differentiability of measures relative
to linear and non-linear operators on $X$. 
Characteristic functionals of measures are studied.
Moreover, the non-Archimedean
analogs of the Bochner-Kolmogorov and Minlos-Sazonov theorems
are investigated. Infinite products of measures are considered and
the analog of the Kakutani theorem is proved.
Convergence of quasi-invariant and pseudo-differentiable
measures in the corresponding spaces of measures
is investigated.
\end{abstract}
\section{Introduction.}
This part is the continuation of the first one and treats
the case of measures with values in non-Archimedean fields of zero
characteristic, for example, the field $\bf Q_p$ of $p$-adic numbers.
There are specific features with formulations of definitions and theorems and
their proofs, because of differences in the notions of $\sigma $-additivity
of real-valued and $\bf Q_p$-valued measures, differences in 
the notions of spaces of integrable functions, quasi-invariance
and pseudo-differentiability.
The Lebesque convergence theorem has quite another meaning, the Radon-Nikodym
theorem in its classical form is not applicable to the considered here case.
A lot of definitions and theorems given below are the non-Archimedean
analogs of classical results. Frequently their formulations and
proofs differ strongly. If proofs differ slightly from the classical
or that of Part I, only general
circumstances are given in the non-Archimedean case.
\par In \S 2 weak distributions, 
characteristic functions of measures and their properties
are defined and investigated.
The non-Archimedean analogs of the Minlos-Sazonov and 
Bochner-Kolmogorov theorems are given. Quasi-measures 
also are considered. 
In \S 3 products of measures are considered together with their 
density functions. The non-Archimedean analog of the Kakutani theorem is
investigated.
In the present paper broad classes of quasi-invariant measures 
are defined and constructed.
Theorems about quasi-invariance of measures under definite 
linear and non-linear transformations
$U: X\to X$ are proved. \S 4 contains a notion of pseudo-differentiability
of measures. This is necessary, because for functions
$f: {\bf K}\to \bf Q_s$ with $s\ne p$ there is not
any notion of differentiability (there is not such non-linear
non-trivial $f$), where $\bf K$ is a field such that ${\bf K}\supset \bf Q_p$. 
There are given criteria for the pseudo-differentiability. 
In \S 5 there are
given theorems about converegence of measures with
taking into account their quasi-invariance and pseudo-differentiability,
that is, in the corresponding spaces of measures.
The main results are Theorems 
2.21, 2.30, 3.5, 3.6, 3.15, 3.19, 3.20, 4.2, 4.3, 4.5, 5.7-5.10.
\par In this part notations of Part I are used also.
\par {\bf Notations.} Henceforth, $\bf K$
denotes a locally compact infinite field with a non-trivial norm,
then the Banach space $X$ is over $\bf K$. 
In the present article measures on $X$ have values 
in the field $\bf K_s$, that is, a finite
algebraic extension of the $s$-adic field $\bf Q_s$ 
with the certain prime number $s$.
Henceforth, $\bf C_s$ denotes the uniform completion of the union
of all $\bf K_s$ with the multiplicative ultranorm
extending that of $\bf Q_s$.
 We assume that $\bf K$ is $s$-free as the additive group,
for example, either $\bf K$ is a finite algebraic extension of 
the field of $p$-adic numbers $\bf Q_p$ or 
$char ({\bf K})=p$ and $\bf K$ is isomorphic with a field
$\bf F_p{(\theta )}$ of formal power series
consisting of elements $x=\sum_ja_j\theta ^j$,
where $a_j\in \bf F_p$, $|\theta |=p^{-1}$, 
$\bf F_p$ is a finite field of $p$ elements,
$p$ is a prime number and $p\ne s$. These imply that 
$\bf K$ has the Haar measures with values in $\bf K_s$
\cite{roo}. If $X$ is a Hausdorff topological space
with a small inductive dimension $ind (X)=0$, then
\par $E$ denotes an algebra of subsets of $X$, 
as a rule $E\supset Bco(X)$ for $\bf K_s$-valued measures, where
\par $Bco (X)$ denotes an algebra of clopen (closed and open)
subsets of $X$,
\par $Bf(X)$ is a Borel $\sigma $-field of $X$ in \S 2.1;
\par $Af(X,\mu )$ is the completion of $E$ by a 
measure $\mu $ in \S 2.1;
\par ${\sf M}(X)$ is a space of norm-bounded measures on $X$ in \S 2.1;
\par ${\sf M}_t(X)$ is a space of Radon norm-bounded measures in \S 2.1;
\par $L(X,\mu ,{\bf K_s})$ is a space of $\mu $-integrable
$\bf K_s$-valued functions on $X$ in \S 2.4;
\par $\chi _{\xi }$ is a character with values in 
$\bf T_s$ in \S 2.5;
\par $\theta (z)=\hat \mu $ is a characteristic functional in \S 2.5;
\par ${\sf C}(Y,\Gamma )$, $\tau (Y)$ in \S 2.20;
\par $\nu \ll \mu $, $\nu \sim \mu $, $\nu \perp \mu $ in \S 2.31.
\section{Weak distributions and families of measures.}
\par  {\bf 2.1.} For a Hausdorff topological space $X$ with
a small inductive dimension $ind (X)=0$ \cite{eng},
henceforth, measures $\mu $ are given on a measurable space $(X,E)$,
where $E$ is an algebra such that
$E\supset Bco(X)$, $Bco(X)$ is an algebra of closed and at the same time
open (clopen) subsets in $X$.
\par We recall that a mapping $\mu : E \to \bf K_s$ 
for an algebra $E$ of subsets of $X$ is called a measure,
if the following conditions are accomplished:
$$(i)\mbox{ }\mu \mbox{ is additive and } \mu (\emptyset )=0,$$
$$(ii) \mbox{ for each }A \in E\mbox{ there exists the following norm}$$
$\| A\|_{\mu }:=\sup \{|\mu (B)|_{\bf K_s}:
\mbox{  }B \subset A, B\in E\}< \infty $,
$$(iii) \mbox{ if there is a shrinking family } F,
\mbox{ that is, for each }$$
$A, B \in F$ there exist
$F\ni C \subset (A\cap B)$ and $\cap\{A: A \in F\}= \emptyset $, then
$\lim_{A\in F} \mu (A)=0$ (see chapter 7 \cite{roo} and also about 
the completion
$Af(X, \mu )$ of the algebra $E$ by the measure $\mu $).
A measure with values in $\bf K_s$ is called a probability measure if
$\| X \|_{\mu } =1$ and $\mu (X)=1$. For functions $f: X\to \bf K_s$ and
$\phi : X\to [0, \infty )$
there are used notations $\| f \|_{\phi }:=\sup_{x \in X} (|f(x)|
\phi (x))$,  $N_{\mu }(x):=\inf(\| U\|_{\mu }: \mbox{ }U \in Bco(X),
\mbox{ } x \in X)$. Tight measures (that is, measures defined on
$E\supset Bco(X)$) compose the Banach space ${\sf M}(X)$ with a norm
$\| \mu \|:=\|X\|_{\mu }$.
Everywhere below
there are considered measures with 
$\|X\|_{\mu }<
\infty $ for $\mu $ with values in $\bf K_s$,
if it is not specified another.
\par A measure $\mu $ on $E$ is called Radon,
if for each $\epsilon >0$ there exists a compact subset $C\subset X$ 
such that $ \| \mu |_{(X\setminus C)} \| < \epsilon $.
Henceforth, ${\sf M}(X)$ denotes the space of norm-bounded
measures, ${\sf M}_t(X)$ is its subspace of Radon
norm-bounded measures.
\par  {\bf 2.2.} If $A \in Bco(L)$, 
then $P_L^{-1}(A)$ is called a cylindrical subset
in $X$ with a base $A$, $B^L:=P_L^{-1}(Bco(L))$,
$B_0:=\cup (B^L:$ $ L\subset X, L  \mbox{ is a Banach subspace },$ 
$dim_{\bf K}X< \aleph _0 )$ (see \S I.2.2).
Let an increasing sequence of Banach subspaces 
$L_n \subset L_{n+1}\subset ...$ such that
$cl(\cup [L_n: n])=X$, $dim_{\bf K}L_n=\kappa _n $ for each $n$
be chosen, where
$cl(A)=\bar A$ denotes a closure of $A$ in $X$ for $A\subset X$. 
We fix a family of projections $P^{L_m}_{L_n}: L_m\to L_n$
such that $P^{L_m}_{L_n}P^{L_n}_{L_k}=P^{L_m}_{L_k}$
for each $m\ge n\ge k$.
A projection
of the measure $\mu $ onto $L$ denoted by
$\mu _L(A):=\mu (P_L^{-1}(A))$ for each $A\in Bco(L)$ compose the consistent
family:
$$(1) \mbox{  } \mu_{L_n}(A)=\mu _{L_m}(P_{L_n}^{-1}(A)\cap L_m)$$
for each 
$m \ge n$, since there are projectors $P_{L_n}^{L_m}$,
where $\kappa _n \le \aleph _0$ and there may be chosen
$\kappa _n< \aleph _0$ for each $n$.
\par An arbitrary family of measures $ \{ \mu _{L_n}: n \in {\bf N} \} $ 
having property $(1)$
is called a sequence of a weak distribution
(see also \cite{dal,sko}).
\par  {\bf 2.3. Lemma.} {\it A sequence of a weak distribution
$\{ \mu _{L_n}: n\}$
is generated by some measure $\mu $ on $Bco(X)$ 
if and only if for each $c>0$ there exists $b>0$ such that
$\| L_n \setminus B(X,0,r) \| _{\mu _{L_n}} \le c$ and
$\sup_n \|L_n\|_{\mu _{L_n}} < \infty $ for $\mu $ with values in
$\bf K_s$, where $r\ge b$.}
\par  {\bf Proof.} For $\mu $ with values in $\bf K_s$ 
the necessity is evident.
To prove the sufficiency it remains only to verify property
(2.1.iii), since then $\|X\|_{\mu }=\sup_n \| L_n\| _{\mu _{L_n}}
< \infty $. Let $B(n) \in E(L_n)$, $A(n)=P_{L_n}^{-1}(B(n))$,
by Theorem 7.6 \cite{roo} for each $c>0$ there is a compact subset
$C(n) \subset B(n)$ such that
$\| B(n)\setminus C(n)\| _{\mu _{L_n}}<c$, where $\| B(n)\setminus
D(n) \|_{\mu } $ $\le \max(\| B(m)\setminus C(m)\| _{\mu _{L(m)}}: m=1,...,n)
< c$ and $D(n):=\bigcap_{m=1}^n {P_{L(m)}}^{-1}(C(m))\cap L_n)$,
$P_{L_n}^{-1}(E(L_n)\subset E=E(X)$. If 
$A(n)\supset A(n+1)\supset ...$ and $\bigcap_n A(n)=\emptyset $, then
$A'(n+1)\subset A'(n)$ and $\bigcap_n A'(n)=\emptyset $, where $A'(n):=
P_{L_n}^{-1}(D(n))$, hence $\| A(n)\| _{\mu } \le \| A'(n)
\| _{\mu }+c$. There
may be taken $B(n)$ as closed subsets in $X$.
In view of the Alaoglu-Bourbaki theorem (see Exer. 9.202(a.3)
\cite{nari}) and the Hahn-Banach theorem (4.8 \cite{roo}) sets
$A(n)$ and $B(X,0,r)$ are weakly compact in $X$, hence, for each $r>0$
there exists $n$ with $B(X,0,r)\cap A(n)=\emptyset $. Therefore,
$\| A(n)\|_{\mu}=\| B(n)\| _{\mu _{L_n}}$ $\le \| L_n\setminus
B(X,0,r)\| _{\mu _{L_n}} \le c$
and there exists $\lim_{n\to \infty } \mu (A(n))=0$,
since $c$ is arbitrary.
\par {\bf 2.4. Definition and notations.} A 
function $\phi :X\to \bf K_s$
of the form $\phi (x)=\phi _S(P_S x)$ is called a cylindrical function if
$\phi _S$ is a $E(S)$-measurable function on a finite-dimensional over
$\bf K$ space $S$ in $X$. For $\phi _S \in L(S,\mu _S,{\bf K_s}):=
L(\mu _S)$ for $\mu $ with values in $\bf K_s$ we may define an integral 
by a sequence of a weak distribution $\{ \mu _{S(n)} \} $: 
$$\int_X \phi (x) \mu _* (dx)
:=\int \phi _{S(n)}(x)\mu _{S(n)} (dx) ,$$ 
where $L(\mu )$ is the Banach space of classes of
$\mu $-integrable functions ($f=g$  $\mu $-almost everywhere, that is,
$\| A\|
_{\mu }=0$, $A:= \{ x:$  $f(x)\ne g(x) \} $ is $\mu $-negligible)
with the following norm $\| f\| :=\| g\| _{N_{\mu }}$ \cite{boui,roo,sko}.
\par {\bf 2.5. Remarks and definitions.} 
In the notation of \S I.2.6
all continuous characters $\chi : {\bf K}\to \bf C_s$
have the form 
$$(1)\mbox{ }\chi _{\xi }(x)=\epsilon ^{z^{-1}\eta ((\xi ,x))}$$ 
for each $ \eta ((\xi ,x)) \ne 0$,
$\chi _{\xi }(x):=1$ for $ \eta ((\xi ,x))=0,$
where $\epsilon =1^z$ is a root of unity, $z=p^{ord(\eta ((\xi ,x)))},$ 
$\pi _j: {\bf K}\to \bf R$,
$\eta (x):= \{ x \} _p$ and $\xi \in {\bf Q_p^n}^*=
\bf Q_p^n$ for $char ({\bf K})=0$,
$\eta (x):=\pi _{-1}(x)/p$ and $\xi \in {\bf K}^*=\bf K$
for $char({\bf K})=p>0$, $x \in \bf K$, 
(see \S 25 \cite{hew}). 
Each $\chi $ is locally constant, hence 
$\chi : {\bf K}\to \bf T_s$ is also continuous, 
where $\bf T$ denotes
the discrete group of all roots of $1$, $\bf T_s$
denotes its subgroup of elements with orders that are not degrees
$s^m$ of $s$, $m \in \bf N$. 
\par For a measure $\mu $ with values in $\bf K_s$
there exists a characteristic functional (that is, called the
Fourier-Stieltjes transformation) $\theta =\theta _{\mu }: C(X,{\bf K}) 
\to \bf C_s$:
$$(2) \mbox{  } \theta (f):=\int _X \chi _e(f(x)) \mu (dx),$$
where either $e=(1,...,1)\in \bf Q_p^n$ for $char ({\bf K})=0$ or 
$e=1\in \bf K^*$ for $char ({\bf K})=p>0$, $x \in X$, $f$ is in the space $C(X,{\bf K})$
of continuous functions from $X$ into $\bf K$, in particular for
$z=f$ in the topologically conjugated space $X^*$
over $\bf K$, $z:X \to \bf K$, $z \in X^*$, $\theta (z)=:
\hat \mu (z)$. 
It has the folowing properties:
$$(3a) \mbox{ }\theta (0)=1 \mbox{ for } \mu (X)=1$$ and $\theta (f)$
is bounded on $C(X,{\bf K})$;
$$(3b)\mbox{ }\sup_f | \theta (f)|=1\mbox{ for probability measures };$$ 
$$(4) \mbox{ } \theta (z) \mbox{ is weakly continuous, that is, } (X^*,
\sigma (X^*,X))\mbox{-continuous},$$
$\sigma (X^*,X)$ denotes a weak topology on $X^*$,
induced by the Banach space $X$ over $\bf K$. 
To each $x \in X$ there corresponds a continuous
linear functional $x^*:X^* \to \bf K$, $x^*(z):=z(x)$, moreover, $\theta
(f) $ is uniformly continuous relative to the norm on 
$$C_b(X,{\bf K}):=\{ f\in C(X,{\bf K}):
\| f\| :=\sup_{x\in X}|f(x)|_{\bf K}<\infty \} .$$
\par  Property (4) follows from Lemma 2.3, boundedness and continuity of
$\chi _e$ and the fact that due to the Hahn-Banach theorem there is
$x_z \in X$
with $z(x_z)=1$ for $z \ne 0$ such that $z|_{(X \ominus L)}=0$ and
$$\theta (z)=\int_X \chi _e(P_L(x)) \mu(dx)=\int_L \chi _e(y) \mu _L(dy),$$
where $L={\bf K}x_z$, also due to the Lebesgue theorem 
(from Exer. 7.F \cite{roo} for $\mu $ with values in
$\bf K_s$). Indeed, for each $c>0$ there exists a
compact subset $S \subset X$ such that $\| X\setminus S\| _{\mu }<c$,
each bounded subset $A \subset X^*$ is uniformly equicontinuous
on $S$ (see (9.5.4) and Exer. 9.202 \cite{nari}), that is,
$ \{ \chi _e(z(x)): $ $z \in A \} $ is the uniformly equicontinuous
family (by $x \in S$). On the other hand, $\chi _e(f(x))$
is uniformly equicontinuous on a bounded $A\subset C_b(X,{\bf K})$ by
$x\in S$.
\par   We call a functional $\theta $ finite-dimensionally concentrated,
if there exists $L \subset X$, $dim_{\bf K}L< \aleph _0 ,$ such that
$\theta |_{(X\setminus L)}=\mu (X)$. For each $c>0$ and $\delta >0$
in view of Theorem 7.6\cite{roo} there 
exists a finite-dimensional over $\bf K$ subspace $L$ and compact
$S \subset L^{\delta }$ such that 
$\| X\setminus S\| _{\mu } <c$. Let
$\theta ^L (z):=\theta (P_Lz)$. 
\par This definition is correct, since
$L\subset X$, $X$ has the isometrical embedding into $X^*$ as
the normed space associated with the fixed basis of $X$, 
such that functionals $z \in X$ separate points in $X$.
If $z \in L$, then $|\theta (z)-\theta ^L(z)| \le c\times b \times q$,
where $b=\| X\| _{\mu }$, $q$ is independent of $c$ and $b$.
Each characteristic functional $\theta ^L(z)$ is uniformly continuous by
$z \in L$ relative to the norm $\| *\|$ on $L$, since $|\theta ^L(z)
-\theta ^L(y)| $ $\le |\int_{S' \cap L} [\chi _e(z(x))-\chi _e(y(x))]$
$\mu _L(dx)|$ $+|\int_{L\setminus S'} [\chi _e(z(x))-\chi _e(y(x))]$
$\mu _L(dx)|$, where the second term does not exceed $2C'$ for 
$ \| L\setminus S'\| _{\mu _L}<c'$ for
a suitable compact subset $S' \subset X$ and $\chi _e(z(x))$ is 
an uniformly equicontinuous
by $x\in S'$ family relative to $z\in B(L,0,1)$.
\par Therefore,
$$(5) \mbox{  } \theta (z)=\lim _{n \to \infty } \theta _n(z)$$
for each finite-dimensional over $\bf K$ subspace $L$, where $\theta
_n(z)$ is uniformly equicontinuous and finite-dimensionally concentrated
on $L_n \subset X$, $z \in X$, $cl(\bigcup_nL_n)=X$, $L_n
\subset L_{n+1}$
for every $n$, for each $c>0$ there are $n$ and $q>0$
such that $|\theta (z)-\theta _j(z)| \le cbq$
for $z \in L_j$ and $j>n$, $q=const >0$ is independent of
$j$, $c$ and
$b$. Let $ \{ e_j:$ $j \in {\bf N} \} $ be the standard orthonormal 
basis in $X$,
$e_j=(0,...,0,1,0,...)$ with $1$ in $j$-th place. Using
Property 2.1.(iii) of $\mu $, local constantness of $\chi _e$,
considering all $z=be_j$ and $b \in \bf K$, we get that $\theta (z)$ on $X$
is non-trivial, whilst $\mu $ is a non-zero measure, since due to Lemma
2.3 $\mu $ is characterized uniquely by $ \{ \mu _{L_n}: n \} $. Indeed,
for $\mu $ with values in $\bf K_s$ a measure $\mu _V$ on $V$,
$dim_{\bf K}V< \aleph _0 ,$ this follows from
Theorem 9.20\cite{hew}, where 
$$F(g)(z):=
\lim _{r \to \infty } \int_{B(V,0,r)} \chi _e(z(x))g(x)m(dx),$$
$z\in V,$ $g \in L(V,\mu _V,{\bf C_s})$, $m$ is the Haar measure on $V$ 
with values in $\bf K_s$.
Therefore, the mapping $\mu \mapsto \theta _{\mu }$ is injective.
\par {\bf 2.6. Theorem.}{\it Let $\mu _1$ and $\mu _2$ be 
measures in ${\sf M}(X)$ on the same algebra
$E$, where $Bco(X)\subset E\subset Bf(X)$
such that $\hat \mu _1(f)=\hat \mu _2(f)$
for each $f \in \Gamma $. Then $\mu _1 =\mu _2$, where $X=c_0(\alpha ,K),$
$\alpha \le \omega _0$, $\Gamma $ is a vector subspace in a space
of continuous functions $f: X\to \bf K$ separating points in $X$.}
\par {\bf Proof.} Let at first $\alpha < \omega _0$, then due to
\S 2.5 $\mu _1=\mu _2,$ since the family $\Gamma $ generates $E$.
Now let $\alpha =\omega _0$, $A=\{x \in
X:\mbox{ } (f_1(x),...,f_n(x)) \in S \},$ $\nu _j$ be an image of a measure
$\mu _j$ for a mapping $x \mapsto (f_1(x),...,f_n(x))$, where 
$S\in E({\bf K^n})$, $ f_j \in
X \hookrightarrow X^*$. Then $\hat \nu _1(y)=\hat \mu _1(y_1f_1+
...+y_nf_n)=\hat \mu _2(y_1f_1+...+y_nf_n)=\hat \nu _2(y)$ for each
$y=(y_1,...,y_n) \in \bf K^n$, consequently, $\nu  _1=\nu _2$ on $E$.
Further compositions of $f\in \Gamma $ with continuous functions
$g: {\bf K}\to \bf K_s$ generate a family of $\bf K_s$-valued 
functions correspondingly separating points of $X$
(see also Chapter 9 in \cite{roo}).
\par  {\bf 2.7. Proposition.} {\it Let $\mu _l$ and $\mu $
be measures in ${\sf M}(X_l)$ and ${\sf M}(X)$ respectively,
where $X_l=
c_0(\alpha _l,K),$ $\alpha _l \le \omega _0,$ $X=\prod_1^n X_l,$
$n \in \bf N.$ Then the condition $\hat \mu (z_1,...,z_n)=\prod_{l=1}^n
\hat \mu _l(z_l)$ for each $(z_1,...,z_n)\in X\hookrightarrow X^*$
is equivalent to $\mu =\prod_{l=1}^n \mu _l$.}
\par {\bf Proof.} Let $\mu =\prod_{l=1}^n \mu _l$, then
$\hat \mu (z_1,...,z_n)=\int _X \chi _e(\sum z_l(x_l))\prod_{l=1}^n
\mu _l (dx_l)$ $=\prod_{l=1}^n \int_{X_l} \chi _e(z_l(x_l))\mu _l(dx_l).$
The reverse statement follows from Theorem 2.6.
\par {\bf 2.8. Proposition.} {\it Let $X$ be a Banach space over $\bf K$;
suppose $\mu $,
$\mu _1$ and $\mu _2$ are probability measures on $X$. Then
the following conditions are equivalent: $\mu $ is the convolution of two 
measures $\mu _j$, $\mu=\mu _1*\mu _2$, and
$\hat \mu (z)=\hat \mu _1
(z)\hat \mu _2(z)$ for each $z \in X$.}
\par {\bf Proof.} Let $\mu =\mu _1 *\mu _2$. This means by the definition
that $\mu $ is the image of the measure $\mu _1 \otimes \mu _2$ for
the mapping $(x_1,x_2) \mapsto x_1 +x_2,$ $x_j \in X,$ consequently,
$\hat \mu (z)=\int_{X\times X} \chi _e(z(x_1+x_2))$ $(\mu _1 \otimes \mu _2)
(d(x_1,x_2))$ $=\prod_{l=1}^2 
\int_X \chi _e(z(x_l))\mu _l(dx_l)$ $=\hat \mu _1
(z)\hat \mu _2(z).$ On the other hand, if $\hat \mu _1\hat \mu _2=\mu ,$
then $\hat \mu =(\mu _1 *\mu _2)^{\wedge }$ and 
due to Theorem 9.20\cite{roo} for measures with values in
$\bf K_s$, we have $\mu =\mu _1*\mu _2.$
\par {\bf 2.9. Corollary.} {\it Let $\nu $ be a probability measure on
$Bf(X)$ and $\mu *\nu =\mu $ for each $\mu $ with values in the same field,
then $\nu =\delta _0$.}
\par {\bf Proof.} If $z_0 \in X \hookrightarrow X^*$
and $\hat \mu (z_0) \ne 0,$  then from $\hat \mu (z_0)\hat \nu (z_0)=
\hat \mu (z_0)$ it follows that $\hat \nu _0(z_0)=1$. From Property
2.6(5) we get that there exists $m \in \bf N$ with $\hat \mu (z)
\ne 0$ for each $z$ with $\| z \| =p^{-m}$, since $\hat \mu (0)=1$.
Then $\hat \nu (z + z_0)=1$, that is, $\hat \nu |_{(B(X,
z_0,p^{-m}))}=1.$ Since $\mu $ are arbitrary we get
$\hat \nu |_X=1$, that is, $\nu =\delta _0$ due to 
\S 2.5.
\par  {\bf 2.10. Corollary.} {\it Let $X$ and $Y$ be Banach spaces
over $\bf K$,
$\mu $ and $\nu $ be probability measures on $X$ and $Y$
respectively, suppose $T: X\to Y$ is a continuous linear operator.
A measure $\nu $ is an image of $\mu $ for $T$ if and only if
$\hat \nu =\hat \mu \circ T^*,$ where $T^*: Y^* \to X^*$ 
is an adjoint operator.}
\par {\bf Proof} follows from \S 2.5 and \S 2.6.
\par  {\bf 2.11. Proposition.} {\it For a completely regular space
$X$ with $ind(X)=0$ the following statements are accomplished:
\par (a) if $(\mu _{\beta })$ is a bounded net of measures 
in ${\sf M}(X)$ that weakly converges to a measure
$\mu $ in ${\sf M}(X)$, then $(
\hat \mu _{\beta } (f))$ converges to $\hat \mu (f)$ for each continuous
$f: X \to \bf K$; if $X$ is separable and metrizable then $(\hat \mu _{\beta })$
converges to $\hat \mu $ uniformly on subsets that are
uniformly equicontinuous in $C(X,{\bf K})$; 
\par (b) if $M$ is a bounded dense family in a ball of the 
space ${\sf M}(X)$ for 
measures in ${\sf M}(X)$, then a family $(\hat \mu: $ $\mu \in M)$
is equicontinuous on a locally $\bf K$-convex space $C(X,{\bf K})$ in a topology
of uniform convergence on compact subsets $S \subset X$.}
\par {\bf Proof.} (a). Functions $\chi _e(f(x))$ 
are continuous and
bounded on $X$, where $\hat \mu (f)=\int_X \chi _e(f(x))\mu (dx)$. Then
(a) follows from the definition of the weak convergence, 
since $sp_{\bf C_s} \{ \chi _e(f(x)):$ $f\in C(X,{\bf K} \} $
is dense in $C(X,{\bf C_s})$.
\par (b). For each $c>0$ there exists a compact subset 
$S \subset X$ such that $\| \mu 
|_{(X\setminus S)} \| <c/4$ for $\bf K_s$-valued measures. 
Therefore, for $\mu \in \sf M$ and $f \in C(X,{\bf K})$ with
$|f(x)|_{\bf K}<c<1$ for $x \in S$ we get $|\mu (X)-{\hat \mu }(f)|=
|\int_X(1-\chi _e(f(x))\mu (dx)|<c/2$
for $\bf K_s$-valued $\mu $, since for $c<1$ and
$x \in S$ we have $\chi _e(f(x))-\chi _e(-f(x))=0$. 
\par {\bf 2.12. Theorem.} {\it Let $X$ be a Banach space over $\bf K$,
$\eta :\Gamma \to \bf C$ be a continuous positive definite function,
$(\mu
_{\beta })$ be a bounded weakly relatively compact net in the
space ${\sf M}_t(X)$ of Radon norm-bounded measures 
and there exists $\lim _{\beta } \hat \mu _{\beta } (f)=
\gamma (f)$ for each $f \in \Gamma $ 
and uniformly on compact subsets of the completion $\tilde \Gamma $,
where $\Gamma \subset C(X,{\bf K})$
is a vector subspace separating points in $X$. Then
$(\mu _{\beta })$ weakly converges to $\mu \in {\sf M}_t(X)$ with
$\hat \mu |_{\Gamma }=\gamma .$}
\par  {\bf Proof} is analogous to the proof of Theorem IV.3.1\cite{vah}
and follows from Theorem 2.6 above and using the non-Archimedean
Lebesgue convergence theorem (see Chapter 7 in \cite{roo}).
\par  {\bf 2.13. Theorem.} {\it (a). A bounded
family of measures in ${\sf M}({\bf K^n})$ 
is weakly relatively compact if and only if a family
$(\hat \mu: $ $\mu \in M)$ is equicontinuous on $\bf K^n$. 
\par (b). If $(\mu _j:$ $j \in {\bf N})$ 
is a bounded sequence of measures
in ${\sf M}_t({\bf K^n})$, $\gamma : {\bf K^n} \to \bf C_s$
is a continuous function, 
$\hat \mu _j(y) \to \gamma (y)$ for each $y \in \bf K^n$
uniformly on compact subsets in $\bf K^n$,
then $(\mu _j)$ weakly converges to a measure
$\mu $ with $\hat \mu =\gamma .$ 
\par (c). A bounded sequence of measures
$(\mu _j)$ in ${\sf M}_t({\bf K^n})$ 
weakly convereges to a measure
$\mu $ in ${\sf M}_t({\bf K^n})$ if and only if for each $y \in \bf K^n$
there exists $\lim _{j \to \infty }\hat \mu _j(y)=\hat \mu (y).$
\par (d). If a bounded net $(\mu _{\beta })$ in ${\sf M}_t({\bf K^n})$
converges uniformly on each bounded subset
in $\bf K^n$, then $(\mu _{\beta })$ converges weakly to a
measure $\mu $ in ${\sf M}_t({\bf K^n})$, where $n \in \bf N$.}
\par {\bf Proof.} (a). This follows from Proposition 2.11. 
\par (b).
Due to the non-Archimedean Fourier transform and the Lebesgue 
convergence theorem \cite{roo} for $\bf K_s$-valued measures and
from the condition $\lim_{R\to \infty } \sup_{|y|>R} | \gamma (y)|R^n=0$
it follows, that for each $\epsilon >0$ there exists $R_0>0$ such that 
$\lim_m \sup_{j>m} \| \mu _j |_{ \{ x\in {\bf K^n}: |x|>R \} } \| 
\le 2\sup_{|y|>R}|\gamma (y)|R<\epsilon $ for each $R>R_0$.
In view of Theorem 2.12 $(\mu _j)$ converges weakly to
$\mu $ with
$\hat \mu =\gamma $. 
(c,d). These can be proved analogously to IV.3.2\cite{vah}.
\par  {\bf 2.14. Corollary.} {\it If $({\hat \mu }_{\beta })\to 1$ 
uniformly on some
neighbourhood of $0$ in $\bf K^n$ for a bounded net
of measures $\mu _{\beta } $ in ${\sf M}_t({\bf K^n})$, then
$(\mu _{\beta })$ converges weakly to $\delta _0$.}
\par  {\bf 2.15. Definition.}  A family of probability measures
$M \subset
{\sf M}_t(X)$ for a Banach space 
$X$ over $\bf K$ is called planely concentrated if
 for each $c>0$ there exists a $\bf K$-linear subspace $S \subset X$ with
$dim_{\bf K}S=n< \aleph _0 $ such that $\inf (\| S^c\| _{\mu }: $ 
$\mu \in M)>1-c$. The Banach space ${\sf M}_t(X)$ is 
supplied with the following norm
$\| \mu \|$
\par {\bf 2.16. Theorem.} {\it Let $X$ be a Banach space 
over $\bf K$ with a family
$\Gamma \subset X^*$ separating points in $M\subset {\sf M}_t(X)$. Then $M$
is weakly relatively compact if and only if a family $ \{ \mu _z:$
$\mu \in M \} $ is weakly relatively compact for each $z \in \Gamma $
and $M$ is planely concentrated, where $\mu _z$ is an image measure 
on $\bf K$ of a measure $\mu $ induced by $z$.}
\par  {\bf Proof} follows from the Alaoglu-Bourbaki theorem \cite{nari},
Lemmas I.2.5 and I.2.21.
\par {\bf 2.17. Theorem.} {\it For $X$ and $\Gamma $ the same as
in Theorem 2.16 a sequence $ \{ \mu _j:$ $j \in {\bf N} \} \subset
{\sf M}_t(X)$  is weakly convergent to $\mu \in {\sf M}_t(X)$
if and only if for each $z \in \Gamma $ there exists
$\lim _{j\to \infty } {\hat \mu }_j(z)={\hat \mu }(z)$ and a family
$\{ \mu _j \} $ is planely concentrated.}
\par  {\bf Proof}  follows from Theorems 2.12,13,16.
\par  {\bf 2.18. Proposition.} {\it Let $X$ be a weakly regular
space with $ind(X)=0$, $\Gamma \subset C(X,{\bf K})$ be a vector subspace
separating points in $X$, $(\mu _n:$ $ n \in {\bf N})$ $\subset
{\sf M}_t(X),$ $\mu \in {\sf M}_t(X)$, $\lim_{n \to \infty }
\hat \mu _n(f)=\hat \mu (f)$ for each $f \in \Gamma $. Then
$(\mu _n)$ is weakly convergent to $\mu $ relative to the weakest
topology $\sigma (X,\Gamma )$ in $X$ relative to which all
$f \in \Gamma $ are continuous.}
\par  {\bf Proof} follows from Theorem 2.13.
\par  {\bf 2.19.} Let $(X,{\sf U})$ $=\prod_{\lambda }
(X_{\lambda },{\sf U}_{\lambda })$ be a product of measurable
completely regular Radon spaces $(X_{\lambda },{\sf U}_{\lambda })$
$=(X_{\lambda },{\sf U}_{\lambda },{\sf K}_{\lambda })$, where
${\sf K}_{\lambda }$ are compact classes approximating from below
each measure $\mu _{\lambda }$ on $(X_{\lambda },{\sf U}_{\lambda })$,
that is, for each
$c>0$  and elements $A$ of an algebra ${\sf U}_{\lambda }$
there is $S \in {\sf K}_{\lambda }$, $S \subset A$
with $\| A\setminus S\| _{\mu _{\lambda }} <c$. 
\par {\bf Theorem.} {\it Each
bounded quasi-measure $\mu $ with values in $\bf K_s$ on $(X,{\sf U})$
(that is, $\mu |_{\sf U_{\lambda }}$ is a bounded measure 
for each $\sf \lambda $)
is extendible to a measure on an algebra $Af(X,\mu ) \supset \sf U$,
where an algebra $\sf U$ is generated by a family $({\sf U}_{\lambda }:$
$\lambda \in\Lambda )$.}
\par {\bf Proof.} We have 2.1(i) by the condition and $\| X\| _{\mu }
< \infty $, if 2.1(iii) is satisfied. It remains to prove 2.1(iii).
For each sequence $(A_n) \subset \sf U$ with $\bigcap _n
A_n=\emptyset $ and each $c>0$ for each $j \in \bf N$ we choose
$K_j \in \sf K$, where the compact class $\sf K$ is generated by
$({\sf K}_{\lambda })$ (see Proposition 1.1.8\cite{dal}), such that
$K_j\subset A_j$ and $\| A_j\setminus K_j\| _{\mu } <c$. Since
$\bigcap _{n=1}^{\infty } K_n $ $\subset \bigcap _n A_n=\emptyset $,
then there exists $l \in \bf N$ with $\bigcap _{n=1}^l K_n=\emptyset $,
hence $A_l=
A_l\setminus \bigcap _{n=1}^l K_n$ $\subset \bigcup _{n=1}^l
(A_n\setminus K_n)$, consequently, $\| A_l\| _{\mu } \le \max_{n=1,...,l}
(\| A_n\setminus K_n\| _{\mu })<c$. It remains to use Theorem 7.8\cite{roo}
about uniqueness of an extension of a measure.
\par {\bf 2.20. Definition.} Let $X$ be a Banach space 
over $\bf K$, then a mapping
$f:X\to \bf C_s$ is called pseudocontinuous, if its restriction $f|_L$
is uniformly continuous for each $\bf K$-linear
subspace  $L\subset X$ with
$dim_{\bf K}L<\aleph _0$. Let $\Gamma $ be a family of mappings $f:Y\to \bf K$
of a set $Y$ into a field $\bf K$. We denote by ${\sf C}(Y,\Gamma )$ an algebra
of subsets of the form $C_{f_1,...,f_n;E}:= \{ x\in X:$ $(f_1(x),...,f_n(x))
\in S \} $, where $S \in Bco({\bf K^n})$, $f_j \in \Gamma $. 
We supply $Y$ with
a topology $\tau (Y)$ which is generated by a base
$(C_{f_1,...,f_n;E}:$ $f_j \in \Gamma ,$
$ E $ $\mbox{is open in}$ $\bf K^n)$.
\par  {\bf 2.21. Theorem. Non-Archimedean analog of the Bochner-Kolmogorov
theorem.} {\it Let $X$ be a Banach space 
over $\bf K$, $X^a$ be its algebraically
adjoint $\bf K$-linear space (that is, of all linear mappings
$f: X\to \bf K$ not necessarily continuous). A mapping $\theta : X^a\to \bf C_s$
is a characteristic functional of a probability measure
$\mu $ with values in $\bf K_s$ and it is defined on $C(X^a,X)$
if and only if $\theta $ satisfies 
Conditions 2.5(3,5) for $(X^a,\tau (X^a)$ and is
pseudocontinuous on $X^a$.}
\par  {\bf Proof.} (I). For $dim_{\bf K}X=card(\alpha )<\aleph _0$
a space $X^a$ is isomorphic with $\bf K^{\alpha }$, hence the statement of this
theorem for a measure $\mu $ with values in $\bf K_s$ follows from Theorem
9.20\cite{roo} and Theorems 2.6 and 2.13 above, since
$\theta (0)=1$ and $|\theta (z)|\le 1$ for each $z$.
\par  (II). Now let $\alpha =\omega _0$. It remains to show that
the conditions imposed on $\theta $ are sufficient, because their necessity
follows from the modification of \S 2.5 (since $X$ has an algebraic embedding
into $X^a$).
The space $X^a$ is isomorphic with ${\bf K}^{\Lambda }$ which is the space
of all $\bf K$-valued functions defined on the Hamel basis $\Lambda $
in $X$. Let $J$ be a family of all non-void
subsets in $\Lambda $. For each $A \in J$ there exists a functional
$\theta _A: {\bf K}^A \to \bf C$ such that $\theta _A(t)=
\theta (\sum_{y \in A}t(y)y)$ for $t \in {\bf K}^A$. From the conditions
imposed on $\theta $ it follows that $\theta _A(0)=1$, $\theta _A$ is
uniformly continuous and bounded on ${\bf K}^A$, 
moreover, due to 2.5(5) for each $c>0$ there are $n$ and $q>0$ such that
for each $j>n$ and $z \in {\bf K}^A$ the following inequality is satisfied:
$$(i)\mbox{  }|\theta _A(z)-\theta _j(z)| \le cbq,$$
moreover, $L_j \supset {\bf K}^A$, $q$ is independent from $j$, $c$ and $b$.
From (I) it follows that on $Bf({\bf K}^A)$ 
there exists a probability measure
$\mu _A$ such that $\hat \mu _A =\theta _A$. The family of measures
$ \{ \mu _A:$ $A \in J \} $ is consistent and bounded, since
$\mu _A=\mu _E \circ
(P_E^A)^{-1}$, if $A \subset E$, where $P_E^A: {\bf K}^E \to
{\bf K}^A$ are the natural projectors. Indeed, 
this is accomplished due to Conditions
(i), 2.5(5) for $X^a$ and due to Theorem 9.20 \cite{roo}.
\par  In view of Theorem 2.19 on
a cylindrical algebra of the space ${\bf K}^{\Lambda }$
there exists the unique measure $\mu $ such that $\mu _A=\mu
\circ (P^A)^{-1}$ for each $A \in J$, where $P^A: {\bf K}^{\Lambda }\to 
{\bf K}^A$
are the natural projectors. From $X^a={\bf K}^{\Lambda }$ it follows that
$\mu $ is defined on $C(X^a,X)$.
For $\mu $ on $C(X^a,X)$ there exists its extension
on $Af(X,\mu )$ such that $Af(X,\mu )\supset Bco(X)$ (see \S 2.1).
\par  {\bf 2.22.} For $f\in L(X,\mu ,{\bf K_s})$ and $\bf K_s$-valued measure
$\mu $ let 
$$\int_Xf(x)\mu _*(dx)=\lim_{n\to \infty }\int_Xg_n(x)\mu _*(dx)$$
for norm-bounded sequence of cylindrical functions
$g_n$ from $L(X,\mu ,{\bf K_s})$
converging to $f$ uniformly on compact subsets of $X$.
Due to the Lebesgue converging theorem this limit exists and does not depend
on a choice of $\{ g_n: n \} $.
\par  {\bf Lemma.} {\it A sequence of a weak distribution
$(\mu _{L_n})$
of probability Radon measures is generated by a $\bf K_s$-valued
probability measure
$\mu $ on $Bco(X)$ of a Banach space 
$X$ over $\bf K$ if and only if there exists
$$(i)\mbox{  }\lim_{|\xi |\to \infty }\int_X G_{\xi }(x)\mu_*(dx)=1,$$
where $\int_X G_{\xi }(x)\mu _*(dx):=S_{\xi }(\{ \mu _{L_n}:n \})$  and \\
$S_{\xi }(\{ \mu _{L_n} \}):=\lim_{n \to \infty }$ $\int_{L_n}
F_n(\gamma _{\xi ,n})(x)$ $\mu _{L_n}(dx),$
$\gamma _{\xi ,n}(y):=\prod_{l=1}^{m(n)}\gamma _{\xi }(y_l)$, \\
$F_n$ is a Fourier transformation by $(y_1,...,y_n)$, $y=(y_j:$ $j \in \bf N)$,
$y_j \in \bf K$,
$\gamma _{\xi }(y):=C(\xi )s^{-2 \min (0,ord_p(y,\xi ))}$,
$C(\xi )\in \bf K_s$,
$\gamma _{\xi }: {\bf K}\to \bf K_s$, $y, \xi \in \bf K$, 
$\nu _{\xi }({\bf K})=1$,$\nu _{\xi }(dy)=\gamma _{\xi }(dy)w(dy)$,
$w: Bco({\bf K})\to \bf K_s$ is the Haar measure; here
$m(n)=dim_{\bf K}L_n< \aleph _0$, $cl(\bigcup_n L_n)=X=c_0(\omega _0,K)$.}
\par {\bf Proof} is quite analogous to that of \S I.2.30 with the substitution
of $|\int_XG_{\xi }(x)\mu _*(dx)-1|<c/2$ for real-valued measures
on $| \| G_{\xi }(x)\| -1|<c/2$ for $\bf K_s$-valued measures.
\par {\bf 2.23. Notes and definitions.}
Suppose $X$ is a locally convex space over a locally compact field
$\bf K$ with non-trivial non-Archimedean valuation
and $X^*$ is a topologically adjoint space.
For a $\bf K_s$-valued measure $\mu $ on $X$ a completion of
a linear space of characteristic functions $ \{ ch_U:
U\in Bco(X) \} $ in $L(X,\mu ,{\bf K_s})$ is denoted by 
$B_{\mu }(X)$.
Then $X$ is called a $KS$-space
if on $X^*$ there exists a topology $\tau $ such that the 
continuity of $f: X^*\to \bf C_s$ with $\| f \|_{C^0} <\infty $
is necessary and sufficient for $f$ to be a characteristic 
functional of a tight measure of the finite norm. Such topology
is called the $K$-Sazonov type topology.
The class of $KS$-spaces contains
all separable locally convex spaces over $\bf K$.
For example, $l^{\infty }(\alpha ,{\bf K})=c_0(\alpha ,{\bf K})^*$.
In particular we also write $c_0({\bf K}):=c_0(\omega _0,{\bf K})$ and
$l^{\infty }({\bf K}):=l^{\infty }(\omega _0,{\bf K})$,
where $\omega _0$ is the first countable ordinal.
\par Let $n_{\bf K}(l^{\infty },c_0)$ denotes the weakest topology on 
$l^{\infty }$ for which all functionals $p_x(y):=\sup_n|x_ny_n|$ are 
continuous, where $x=\sum_nx_ne_n\in c_0$ and $y=\sum_ny_ne^*_n
\in l^{\infty }$, $e_n$ is the standard base in $c_0$.
Such topology $n_{\bf K}(l^{\infty },c_0)$ is called the normal topology.
The induced topology on $c_0$ is denoted by $n_{\bf K}(c_0,c_0)$.
\par {\bf 2.24. Theorem.} {\it Let $f: l^{\infty }({\bf K})\to \bf C_s$
be a functional such that
\par $(i)$ $f(0)=1$ and $ \| f \| _{C^0}\le 1$,
\par $(ii)$ $f$ is continuous in the normal topology
$n_{\bf K}(l^{\infty },c_0)$, then $f$ is the characteristic functional
of a probability measure on $c_0({\bf K})$.}
\par {\bf Proof.} If $\nu $ is the Haar measure on $\bf K^n$, then
on $Bco({\bf K^n})$ it takes values in $\bf Q$. Therefore, Lemma 4.1
\cite{mamaz} is transferable onto the case of $\bf K_s$-valued measures,
since ${\bf Q}\subset \bf K_s$. Therefore, analogously to Equation
$(4.1)$ of Lemma 4.2 \cite{mamaz} we have
$$(i)\mbox{ }P \{ |V_1|_{\bf K}<\epsilon ,...,
|V_n|_{\bf K}<\epsilon \} = \nu ^{-1}(B({\bf K^n},0,p^{-m}))
\int_{\bf K^n} f_V(y)ch_{B({\bf K^n},0,p^{-m})}(y)\nu (dy)$$
for measurable maps $V_j: (\Omega ,{\sf B},P)\to ({\bf K},Bco({\bf K}))$,
where $(\Omega ,{\sf B},P)$ is a probability space for a probability measure
$P$ with values in $\bf K_s$ on an algebra ${\sf B}$ of subsets of a set 
$\Omega $, $f_W$ is a characteristic function of $W=(V_1,...,V_n)$. 
To continue the proof we need the following statements.
\par {\bf 2.25. Lemma} {\it Let $f: c_0({\bf K})\to \bf C_s$ be a function
satisfying the following two conditions:
\par $(i)$ $|f(x)|\le 1$ for each $x\in c_0({\bf K})$,
\par $(ii)$ $f$ is continuous at zero in the topology 
$n_{\bf K}(c_0,c_0)$,  \\
then for each $\epsilon >0$ there exists $\lambda (\epsilon )\in c_0
({\bf K})$ such that $|1-f(x)|<p_{\lambda (\epsilon )}(x)+\epsilon $
for each $x\in c_0({\bf K}).$}
\par {\bf Proof.} In view of continuity for each $\epsilon >0$
there exists $y(\epsilon )\in c_0$ such that $|1-f(x)|<\epsilon $ 
if $p_{y(\epsilon )}<1$. Put $\lambda (\epsilon )=\pi _{\bf K}^{-1}
y(\epsilon )$, where $\pi _{\bf K}\in \bf K$ is such that 
$|\pi _{\bf K}|=p^{-1}$. If $x\in c_0$ is such that
$p_{\lambda (\epsilon )}(x)<p^{-1}$, then
$|1-f(x)|<\epsilon \le \epsilon + p_{\lambda (\epsilon )}(x).$
If $p_{\lambda (\epsilon )}(x)\ge p$, then $|1-f(x)|\le 2 \le p
<p_{\lambda (\epsilon )}(x)+\epsilon $.
\par {\bf 2.26. Lemma.} {\it Let $\{ V_n : n\in {\bf N} \} $ 
be a sequence of $\bf K$-valued random variables
for $P$ with values in $\bf K_s$.
If for each $\beta >0$ and $\epsilon >0$ there exists 
$N_{\epsilon }\in \bf N$ such that
$$(i)\mbox{ } \| P |_{
\{ \sup_{n\ge N_{\epsilon }} |V_n|_{\bf K}\le \beta  \} } \|
\ge 1-\epsilon (1+\beta ^{-1}),$$ then
$\lim_nV_n=0$ $P$-a.e. on $\Omega $.}
\par {\bf Proof} is quite analogous to that of Lemma 4.4 \cite{mamaz}
with substitution of $P$ on $\| P \| $.
\par {\bf 2.27. Proposition.} {\it Let $f: c_0({\bf K})\to \bf C_s$ 
be a function such that
\par $(i)$ $f(0)=1$ and $|f(x)|\le 1$ for each $x\in c_0$,
\par $(ii)$ $f(x)$ is continuous in the normal topology $n_{\bf K}(c_0,c_0)$.
Then there exists a probability measure $\mu $ on $c_0({\bf K})$
such that $f(x)={\hat \mu }(x)$ for each $x\in c_0$.}
\par {\bf Proof.} Consider functions $f_n(x_1,...,x_n):=f(
x_1e_1+...+x_ne_n),$ where $x=\sum_jx_je_j\in c_0$.
From Condition $(ii)$ and Proposition 3.1(2) \cite{mamaz}
it follows, that $f(x)$ is continuous in the norm topology.
From Chapters 7,9 \cite{roo} it follows, that there
exists a consistent family of tight measures $\mu _n$ on $\bf K^n$ such that
${\hat \mu }_n(x)=f_n(x)$ for each $x\in \bf K^n$.
In view of Theorem 2.19 there exists a probability space
$(\Omega ,{\sf B},P)$ with a $\bf K_s$-valued measure $P$ and a sequence 
of random variables $\{ V_n \} $ such that
$\mu _n(A)=P\{ \omega \in \Omega :$ 
$(V_1(\omega ),...,V_n(\omega ))\in A \} $ for each clopen subset 
$A$ in $\bf K^n$, consequently, $\lim_nV_n=0$ $P$-a.e. in $\Omega $.
In view of the preceding lemmas we have the following inequality:
$$|1- \| P|_{(|V_n|<\beta ,...,|V_{n+m}|<\beta )} \|
\le \| p_{\lambda (\epsilon )} (y_1e_n+...+y_me_{n+m}\| _{
L(B({\bf K^n},0,\beta ^{-1}),\nu ,{\bf K_s})}.$$
Since $\lim_kp_{\lambda (\epsilon )}(e_k)=0$, then there exists
$N\in \bf N$ such that $\sup_{k\ge N}p_{\lambda (\epsilon )}
(e_k)\le \epsilon $,
consequently, $ \| P |_{ \{ |V_N|<\beta ,...,|V_{N+m}||<\beta  \} } \|
\ge 1-\epsilon (1+\beta ^{-1})$. Due to Lemma 2.34
$ \| P |_{ \{ \lim_nV_n=0 \} } \| =1$.
Define a measurable mapping
$W$ from $\Omega $ into $c_0$ by the following formula:
$W(\omega ):=\sum_nV_n(\omega )e_n$ for each $\omega \in \Omega $,
then we also define a measure $\mu (A):=P \{ W^{-1}(B) \} $
for each $A\in Bco(X)$, hence $\mu $ is a probability measure on $c_0$.
In view of the Lebesgue convergence theorem (see Chapter 7 \cite{roo})
there exists ${\hat \mu }(x)=\lim_n{\hat \mu }_n(x_1e_1+...+x_ne_n)=f(x)$
for each $x\in c_0$.
\par {\bf Continuation of the proof of Theorem 2.24.}
Let $f: l^{\infty }({\bf K})\to \bf C_s$ satisfies
assumption of Theorem 2.24,  then by Proposition 2.27
there exists a probability measure $\mu $ on $c_0({\bf K})$
such that $f(x)={\hat \mu }(x)$ for each $x\in c_0({\bf K}).$
\par {\bf 2.28. Theorem.} {\it Let $\mu $ be a probability measure
on $c_0({\bf K})$, then $\hat \mu $ is continuous in 
the normal topology $n_{\bf K}(l^{\infty },c_0)$ on $l^{\infty }$.}
\par {\bf Proof.} It is quite analogous to that of I.2.33
due to Lemma 2.3 and Theorem 2.19.
\par {\bf 2.29. Corollary.} {\it The normal topology $n_{\bf K}
(l^{\infty },c_0)$ is the $K$-Sazonov type topology
on $l^{\infty }({\bf K})$.}
\par  {\bf 2.30. Theorem. Non-Archimedean analog of the
Minlos-Sazonov theorem.}
{\it For a separable Banach space 
$X$ over $\bf K$ the following two conditions are equivalent:
$$(I)\mbox{ }\theta : X \to {\bf C_s}
\mbox{ satisfies Conditions }2.5(3,5)
\mbox{ and }$$ for each $c>0$ there exists a compact operator
$S_c: X\to X$ such that 
$|\theta (y)-\theta (x)|<c$ for $|\tilde z (S_cz)|<1$;
$$(II)\mbox{ }\theta \mbox{ is a characteristic functional of a 
probability Radon measure } \mu $$ on $E$, where $\tilde z$ is an element
$z \in X
\hookrightarrow X^*$ considered as an element of $X^*$ under the
natural embedding associated with the standard base of 
$c_0(\omega _0,{\bf K})$, $z=x-y$, $x$
and $y$ are arbitrary elements of $X$.}
\par {\bf Proof.} $(II\to I)$.
For $\theta $ generated by a $\bf K_s$-valued measure for each $r>0$
we have $| \theta (0)-\theta (x)|=|\int_X(1-\chi _e(x(u)))
\mu (du)|\le \| (1-\chi _e(x(u)))|_{B(X,0,r)} \| _{\mu }+2 \| \mu |_{(X
\setminus B(X,0,r))} \| $. In view of the
Radon property of the space $X$ and Lemma I.2.5 for each $b>0$ and
$\delta >0$ there are a finite-dimensional over $\bf K$ subspace
$L$ in $X$ and a compact subset $W \subset X$ such that
$W \subset L^{\delta }$, $\| \mu |_{(X\setminus W)} \| <b$, 
hence $\| \mu |_{(X\setminus L^{\delta })} \| <b.$
\par We consider the expression $J(j,l)$ (see \S I.2.35).
and the compact operator $S: X \to X$ with $\tilde e_j(Se_l)=
\xi _{j,l}t.$
Then $|\theta (0)-\theta (z)|<c/2 + 
|{\tilde z}(Sz)|<c$ for the $\bf K_s$-valued measure,
if $|\tilde z(Sz)| <|t|c/2$. We choose $r$ such that
$\| \mu |_{(X\setminus B(X,0,r))} \| <c/2$ 
with $S$ corresponding to $(r_j:$ $j)$,
where $r_1=r$,
$L_1=L$, then we take $t \in \bf K$ with $|t|c=2$.
\par  $(I \to II)$. Without restriction of generality we may take
$\theta (0)=1$
after renormalization of non-trivial $\theta $. 
In view of Theorem 2.24 as in \S 2.5 
we construct using $\theta (z)$ a consistent family
of finite-dimensional distributions $\{ \mu _{L_n}: n \} $
all with values in $\bf K_s$. Let $m_{L_n}$
be the $\bf K_s$-valued Haar measure 
on $L_n$ which is considered as $\bf Q_p^a$ with
$a=dim_{\bf K}L_ndim_{\bf Q_p}{\bf K}$, $m(B(L_n,0,1))=1$.
If $S_c$ is a compact operator such that $|\theta (y)-\theta (x)|<c$
for $|\tilde z(S_cz)|<1$, $z=x-y$, then 
$|1-\theta (x)|< \max (C, 2| {\tilde x}(S_cx)|)$ and
$ \| \gamma _{\xi ,n}(z)(1-\theta (z)) \| _{m_{L_n}} \le $ \\ 
$\max ( \| \gamma _{\xi ,n}(z) \| _{m_{L_n}}C,
2|(\gamma _{\xi ,n}(z)){\tilde z}(S_cz)|_{m_{L_n}})\le $
$\max (C,b \| S_c \| /|\xi |^2)$, \\
where $b:=p\times \sup _{|\xi |>r}
(|\xi |^2 \| \gamma _{\xi ,n}(z) z^2 \| _{m_{L_n}})<\infty $ for 
the $\bf K_s$-valued measures.
Due to the formula of changing variables in integrals (A.7\cite{sch1})
the following equality is valid: 
$$|1- \| G_{\xi }(x) \| _{\mu _*}|\le \max (C,b \| S_c \|/|\xi |^2)$$
for the $\bf K_s$-valued measures.
Then taking the limit with $|\xi | \to \infty $ and then
with $c \to +0$ with the help of Lemma 2.22 we get the statement $(I\to II)$.
\par  {\bf 2.31. Definition.} Let on a completely regular space
$X$ with $ind (X)=0$ two non-zero $\bf K_s$-valued
measures $\mu $ and $\nu $ are given.
Then $\nu $ is called absolutely continuous relative to $\mu $ if 
there exists $f$ such that $\nu (A)=\int _A
f(x)\mu (dx)$ for each $A\in Bco(X)$, where $f \in L(X,\mu ,{\bf K_s})$ 
and it is denoted $\nu \ll \mu $. Measures $\nu $ and $\mu $ are
singular to each other if there is $F \in E$ with
$ \| X\setminus F\| _{\mu }=0$ and $\| F\| _{\nu }=0$
and it is denoted
$\nu \perp \mu $. If $\nu \ll \mu $ and $\mu \ll \nu $ then they are called
equivalent, $\nu \sim \mu $.
\par  {\bf 2.32. Definition and note.} For $\mu : E(X) \to \bf K_s$ a sequence
$(\phi _n(x): n) \subset L(\mu )$ is called a martingale if for each
$\psi \in L(\mu | {\sf U_n})$:
$$(i)\mbox{ }\int_X\phi _{n+1}(x)\psi (x) \mu (dx)\mbox{ } =\int_X\phi _n(x)
\psi (x) \mu (dx)$$ such that $(\phi _n:$ $ n)$ is uniformly converging
on $Af(X, \mu )$-compact subsets
in $X$, where $\sf U_n$ is the minimal algebra
such that $(\phi _j :$ $ j=1,...,n)\subset
L(\mu |{\sf U_n})$, $\mu | \sf U_n$ is a restriction of $\mu $ on
${\sf U_n}
\subset E(X)$, $X$ is the Banach space over $\bf K$.
\par In view of \S \S 7.10 and 7.12\cite{roo}
for $\| X\| _{\mu } < \infty $
the $Af(X,\mu )$-topology on compact subspaces $X_c:=[x \in X:$
$N_{\mu }(x) \ge c]$ coincides with the initial topology, if $\mu $
is defined on $E$ such that $Bco(X)\subset E \subset Af(X, \mu )$, where $c>0$.
\par  {\bf 2.33. Theorem.} {\it If there is a martingale $(\phi _n:$ $ n)$
for $\mu $ with values in $\bf K_s$ and $\sup _n \| \phi _n\| _{N_{\mu }}
< \infty $, then there exists $\lim _{n \to \infty } \phi _n(x)=:
\phi (x) \in L(\mu )$.}
\par  {\bf Proof.}  Let $\psi (x)$ be a characteristic function of a
clopen subset in $X$,
then for each $\phi _n$ there exists a sequence of simple functions
$(\phi _n^j:$
$j \in \bf N$ $)$ such that $\lim _{j \to \infty }
\| \phi _n - \phi _n^j\| _{N_{\mu}}=0$. From
$\| \phi _n -\phi _n^{j(n)}\|
_{N_{\mu }}<c$ and 2.32.(i) 
it follows that $|\int_X(\phi _{n+1}^{j(n+1)}(x)-
\phi _n^{j(n)})\psi (x) \mu (dx)| <c\| \psi \| _{N_{\mu }}$ for each
$\psi \in L(\mu )$, consequently, $\| \phi _{n+1}^{j(n+1)}- \phi _n^{j(n)}\|
_{N_{mu }}<c$ and there exists $\lim_{n \to \infty } \phi _n^{j(n)}=
\lim_{n \to \infty } \phi _n=\phi \in L(\mu )$ due to the Lebesgue
theorem, if $(c=c(n)=s^{-n}:$ $n \in {\bf N})$, where for each $\phi _n$
is chosen $j(n) \in \bf N$, since $(\phi _n ^{j(n)}:$ $n)$
is a Cauchy sequence in the Banach space 
$L(\mu )$ due to the ultrametric inequality.
\section{ Quasi-invariant measures.}
\par In this section after few preliminary statements there are given
the definition of a quasi-invariant measure 
and the theorems about quasi-invariance of measures
relative to transformations of a Banach space $X$ over $\bf K$.
\par {\bf 3.1.} Let $X$ be a Banach space over $\bf K$, $(L_n:n)$
be a sequence of subspaces, $cl (\bigcup _n L_n)=X$, $L_n
\subset L_{n+1}$ for each $n$, $\mu ^j$ be probability measures,
$\mu ^2 \ll \mu ^1,$ $(\mu ^j_{L_n})$ be sequences of weak distributions,
also let there exist derivatives
$\rho _n(x)=\mu ^2_{L_n}(dx)/ \mu ^1_{L_n}(dx)$ and the following limit
$\rho (x):=\lim _{n \to \infty }\rho _n(x)$ exists.
\par {\bf Theorem.} {\it If $\mu ^j $ are $\bf K_s$-valued
and $[\rho _n(P_{L_n}x):$ $n]$ converges uniformly on
$Af(X, \mu ^1)$-compact subsets in $X$, $\sup_n\| \rho
_n\| _{N_{\mu ^1}} < \infty $, then this is equivalent to the following:
$\rho (x)=\mu ^2(dx)/ \mu ^1(dx) \in L(\mu ^1)$ and $\lim_{n
\to \infty }\|\rho (x)-\rho _n(P_{L_n}x)\|_{N_{\mu ^1}}=0$.}
\par  {\bf Proof}  For each $A \in Bco(L)$ the equality is accomplished:
$$\mu ^2_L(A)=\int_A \rho _L(x)\mu ^1_L(dx)=\int_{P_{L}^{-1}(A)}
\rho _L(P_Lx) \mu ^1(dx).$$
 Then for each $\psi \in L(\mu ^1|P_L^{-1}[Bco
(L)])$ we have  
$$\int_X\psi(x)\mu ^2(dx)=\int_X\rho _L(P_Lx)\psi (x)
\mu ^1(dx),\mbox{ consequently,}$$
$$\int_X \rho _{n+1}(x)\psi (x)\mu ^1(dx)=\int_X\psi (x)\mu ^2(dx)
=\int_X \rho _n(x)\psi (x) \mu ^1(dx),$$
where $\rho _{L_n}=\rho _n$, $\psi \in L(\mu ^1|P_{L_{n+1}}^{-1}[Bco(L
_{n+1})])$. From Theorem 2.33 and Definition 2.31 the statement follows.
\par  {\bf 3.2. Theorem.} {\it (A). Measures $\mu ^j: E\to \bf K_s$, 
$j=1,2$, for a Banach space $X$ over $\bf K$ are orthogonal
$\mu ^1 \perp \mu ^2$
if and only if 
$N_{{\mu ^1}}(x)N_{{\mu ^2}}(x)=0 $ for each $x \in X$.
\par (B). If for measures $\mu ^j: E \to \bf K_s$ on a
Banach space $X$ over $\bf K$
is satisfied $\rho (x)=0$ for each $x$ with $N_{{\mu ^1}}(x)>0$,
then $\mu ^1\perp \mu ^2$; the same is true for a completely regular
space $X$ with $ind(X)=0$
and $\rho (x)=\mu ^2(dx)/ \mu ^1(dx)=0$ for 
each $x$ with $N_{{\mu ^1}}(x)>0$.}
\par  {\bf Proof.} (A). 
From Definition 2.31 it follows
that there exists $F \in E$ with $\| X\setminus F\|_{{\mu ^1}}=0$ and
$\|F\|_{{\mu ^2}}=0$. In view of Theorems 7.6 and 7.20\cite{sko}
the characteristic function $ch_F$ of the set $F$ belongs to
$L(\mu ^1)\cap L(\mu ^2)$
such that $N_{{\mu ^j}}(x)$ are semi-continuous from above,
$\| ch_F \|_{N_{\mu ^2}}
=0$, $\| ch_{X\setminus F} \|_{N_{{\mu ^1}}}=0$, consequently,
$N_{{\mu ^1}}(x)N_{{\mu ^2(x)}}=0$ for each $x \in X$.
\par  On the other hand, if $N_{{\mu ^1}}(x)N_{{\mu ^2}}(x)=0$ for each
$x$, then for $F:=[x \in X:$ $N_{\mu ^2}(x)=0]$ due to Theorem 7.2 \cite{roo}
$\|F\|_{{\mu ^2}} =\|ch_F\|_{N_{{\mu ^2}}}=0$. 
Moreover, in view of Theorem 7.6\cite{roo} $F=
\bigcap_{n=1}^{\infty } U_{s^{-n}}$, where $U_c:=[x \in X:$ $N_{\mu ^2}
(x)<c]$ are open in $X$, hence $ch_F \in L(\mu ^1)\cap L(\mu ^2)$
and $N_{\mu ^1}|_{(X\setminus F)}=0$, consequently, $\|X\setminus
F\|_{\mu ^1}=0$.
\par  (B). In view of Theorem 2.19 for each $A \in P_{L_n}^{-1}
[E(L_n)]$ and $m>n$: $\int_A \rho _m(x)\mu ^1(dx)=\mu ^2(A)$, then from
$\lim_{n \to \infty }\| \rho (x)- \rho _n(P_{L_n}x)\|
_{N_{\mu ^1}}=0$ and Conditions 2.1.(i-iii) on $\mu ^2$
Statement (B) follows.
\par  {\bf 3.3. Note.}
The Radon-Nikodym theorem is not valid for  $\mu ^j$ with values in
$\bf K_s$, so not all theorems for real-valued measures may be
transferred onto this case. Therefore, the definition
of absolute continuity of measures was changed (see \S 
2.31 and \cite{schrn}).
\par  {\bf 3.4. Theorem.} {\it Let measures $\mu ^j$ and $\nu ^j$
be with values in $\bf K_s$ on $Bco(X_j)$ for a 
Banach space $X_j$ over $\bf K$
and $\mu =
\mu ^1 \otimes \mu ^2$, $\nu =\nu ^1 \otimes \nu ^2$ on $X=X_1
\otimes X_2$, therefore, the statement $\nu \ll \mu $ is equivalent to
$\nu ^1 \ll \mu ^1$ and $\nu ^2 \ll \mu ^2$, moreover,
$\nu (dx)/ \mu (dx)=
(\nu ^1(P_1dx)/ \mu ^1 (P_1dx))(\nu ^2(P_2dx)/ \mu ^2(P_2dx))$, where
$P_j:X \to X_j$ are projectors.}
\par  {\bf Proof} follows from Theorem 7.15\cite{roo}
and modification of the proof of Theorem 5 \S 15\cite{sko}.
\par   {\bf 3.5. Theorem. The non-Archimedean analog of 
the Kakutani theorem.} {\it Let $X=\prod _{j=1}^{\infty } X_j$
be a product of completely regular spaces $X_j$ with
$ind(X_j)=0$ and probability measures $\mu ^j,$ $\nu ^j: E(X_j)\to \bf K_s$,
also let $\mu _j \ll \nu _j$ for each
$j$, $\nu =\bigotimes_{j=1}^{\infty }\nu _j$, $\mu = \bigotimes _{j=1}
^{\infty }\mu _j$ are measures on $E(X)$,
$\rho _j(x)=\mu _j(dx)/\nu _j(dx)$
are continuous by $x \in X_j$, $\prod_{j=1}^n \rho _j(x_j)=:t_n(x)$
converges uniformly on
$Af(X, \mu )$-compact subsets in $X$, $\beta _j
:=\| \rho _j(x)\|_{\phi _j}$, $\phi _j(x):=N_{\nu ^j}(x)$ on $X_j$.
If $\prod _{j=1}^{\infty } \beta _j$ converges in $(0, \infty )$
(or diverges to 0), then $\mu \ll \nu $ and $q_n(x)=\prod _{j=1}^n\rho _j
(x_j)$ converges in $L(X,\nu ,{\bf K_s})$ 
to $q(x)=\prod _{j=1}^{\infty }\rho _j(x_j)$
$=\mu (dx)/ \nu (dx)$ (or $\mu \perp \nu$ respectively), where $x_j
\in X_j$, $x \in X$.}
\par  {\bf Proof.} The countable additivity of $\nu $ and $\mu $
follows from Theorem 2.19. Then $\beta _j=\| \rho _j\| _{\phi _j}
\le \| \rho _j\| _{N_{\nu _j}}=\| X \| _{\mu _j}=1$, since $N_{\nu _j}
\le 1$ for each $x \in X_j$, hence $\prod _{j=1}^{\infty } \beta _j$
can not be divergent to $\infty $. If this product diverges to $0$
then there exists a sequence $\epsilon _b:=\prod_{j=n(b)}^{m(b)}
\beta _j$ for which the series converges $\sum _{b=1}^{\infty } \epsilon _b
<\infty $,
where $n(b) \le m(b)$. For $A_b:=[x:$ $(\prod_{j=n(b)}^{m(b)} \rho _j(x_j)
) \ge 1]$ there are estimates $\| A_b \| _{\nu } \le $ $\sup_{x \in A_b}
[\prod_{j=n(b)}^{m(b)}| \rho _j(x_j)| \phi _j(x_j)] \le \epsilon _b$,
consequently, $\| A \| _{\nu }=0$ for $A=\lim \sup (A_b:$ $b \to \infty )$,
since $0< \sum_{b=1}^{\infty } \epsilon _b < \infty $.
\par  For $B_b:=X\setminus A_b$ we have: $\| B_b \| _{\mu } \le $ $[\sup_{x
\in B_b} \prod _{j=n(b)}^{m(b)} |1/ \rho _j(x_j)| \psi (x_j)]$
$=[\prod_{j=n(b)}^{m(b)} \| \rho _j(x_j) \| _{\phi _j}]=
\epsilon _b$, where $\psi _j(x)=N_{\mu _j}(x)$, since $\mu _j(dx_j)
= \rho _j(x_j) \nu _j(dx_j)$ and $N_{\mu _j}(x)=|\rho _j(x_j)|
N_{\nu _j}(x)$ due to continuity of $\rho _j (x_j)$ (for $\rho _j
(x_j)=0$ we set $|1/ \rho _j(x_j)| \psi _j(x_j)=0$, because
$\psi _j(x_j)=0$ for such $x_j$), consequently, $\| \lim \sup (B_b:$
$ b \to \infty )\| _{\mu }=0$ and $\| A \| _{\mu } \ge \| \lim \inf
(A_b:$ $b \to \infty) \| _{\mu }=1$. This means that $\mu \perp \nu $.
\par  Suppose that $\prod_{j=1}^{\infty } \beta _j$ converges to
$0<\beta < \infty $, then $\beta \le 1$ (see above). 
Therefore from the Lebesgue
Theorem 7.F\cite{roo} it follows that $t_n(x)$ converges in 
$L(X,\mu ,{\bf K_s})$, since
$|t_n(x)| \le 1$ for each $x$ and $n$, at the same time
each $t_n(x)$ converges uniformly on compact subsets in the topology
generated by $Af(X, \mu )$.
Then for each bounded continuous cylindrical function
$f: X\to \bf K_s$ we have 
$$\int_X f(x) \mu (dx)=\int_X f(x_1,...,x_n)
t_n(x) \otimes _{j=1}^n \nu _j(dx_j)=$$ 
$$\lim _{n \to \infty }
\int_X f(x)t_n(x)\nu (dx)=\int_X \rho (x) \nu (dx).$$ Approximating
arbitrary $h \in L(X,\mu ,{\bf K_s})$ by such $f$ we get the equality
$$\int_X h(x) \mu (dx)=\int_X h(x) \rho (x) \nu (dx),$$ consequently,
$\rho (x)=\mu (dx)/\nu (dx)$.
\par  {\bf 3.6. Theorem.} {\it Let $\nu $, $\mu $, $\nu _j$, $\mu _j$
be probability measures with values in $\bf K_s$, $X$
and $X_j$
be the same as in \S 3.5 and $\mu \ll \nu $, then $\mu _j \ll \nu _j$
for each $j$ and $\prod _{j=1}^{\infty } \beta _j$ converges to $\beta $,
$\infty
> \beta > 0$, where 
$\beta _j=\| \rho _j \| _{\phi _j},$ $\phi _j(x)= N_{\nu _j}(x)$.}
\par {\bf Proof.} For $\bf K_s$-valued measures
from $P_j^{-1}(Bco(X_j)) \subset Bco(X)$
it follows that $\mu _j \ll \nu _j$ for each $j$, since $\prod _1^{\infty }
\rho _j(x_j)=\rho (x)$ $\in L(X, \nu )$ and $ \rho _j(x_j) \in L(X_j,
\nu _j)$, where $x_j=P_jx$, $P_j: X \to X_j$ are projectors. Then
$\rho (x)=\lim_{n \to \infty } \prod _1^n\rho _j(P_jx)$ and $\| \rho (x)
\| _{N_{\nu }}=\lim_{n \to \infty } \| \rho _j \| _{N_{\nu _j}}$.
Since $N_{\nu _j} \le 1$, then $\phi _j(x) \le N_{\nu _j}(x)$ and
for $\phi =N_{\nu }$, consequently, $\| \rho (x) \| _{\phi }$
$=\lim _{n \to \infty } \prod _{j=1}^n \| \rho _j \| _{\phi _j}$
$\le \| \rho \| _{N_{\nu }} =1$ (due to the definition of the Tihonov
topology in $X$ [see \S 2.3\cite{eng}] and definition of $\| *\| _{\phi }$).
If $\| \rho \| _{\phi }=0$, then $\| \rho \| _{N_{\nu }}=0$
and by Theorem 3.2(B) this would mean that $ \nu \perp \mu $ or
$\mu =0$, but $\mu \ne 0$, hence $\beta >0$.
\par {\bf 3.7. Definition.}
Let $X$ be a Banach space over $\bf K$, $Y$ be a completely regular space with
$ind(X)=0$, $\nu : Bco(Y) \to \bf K_s$,
$\mu ^y: Bco(X) \to \bf K_s$ for each $y \in Y$, suppose $\mu ^y(A)
\in L(Y, \nu )$ for each $A \in Bco(X)$, $\| Y \| _{\nu }< \infty $,
$\sup _{y \in Y} \| X \| _{\mu ^y}< \infty $, a family
$( \mu ^y(A_n):$ $n)$ is converging uniformly by $y\in C$ on each
$Af(Y, \nu )$-compact subset $C$ in $Y$
for each given shrinking family of subsets $(A_n:$ $n)\subset X$.
Then we define:
$$(i)\mbox{ } \mu (A)=\int_Y \mu _y(A) \nu (dy).$$
A measure $\mu $ is called mixed. 
Evidently, Condition 2.1(i) is fulfilled;
(ii): $\| A \| _{\mu } \le
(\sup_{y \in Y} \| A \| _{\mu ^y})\| A \| _{\nu } < \infty $;
(iii) is carried out due to the Lebesgue theorem, since
$\lim_{n \to \infty }
\mu (A_n)$ $=\int_Y(\lim_n \mu ^y(A_n))\nu (dy)=0$.
We define measures $\pi ^j$ by the formula:
$$(ii)\mbox{ }\pi ^j(A\times C)=\int_C \mu ^{j,y}(A) \nu ^j(dy),$$
where $j=1,2$ and $\mu ^{y,j}$ together with $\nu ^j$ are defined
as above $\mu ^y$ and $\nu $. 
\par  {\bf 3.8. Theorem.} {\it Let $\mu ^j$ be $\bf K_s$-valued
measures and $\pi ^j$, $X$ and $Y$
be the same as in \S 3.7.
\par (A). If $\pi ^2 \ll \pi ^1$, then $\nu ^2
\ll \nu ^1$ and $\mu ^{2,y} \ll \mu ^{1,y}$ $(mod$ $\nu ^2)$.
\par  (B). If $\nu ^2 \ll \nu ^1$ and $\mu ^{2,y} \ll \mu ^{1,y}$ $(mod$
$\nu ^2)$ and a $Bco(X\times Y, \pi ^1)$-measurable function
$\tilde \rho (y,x) =\mu ^{2,y}(dx)/ \mu ^{1,y}(dx) \in
L(X\times Y, \pi ^1)$ exists, then $\pi ^2 \ll \pi ^1$ and
$\pi ^2(d(x,y))/\pi ^1(d(x,y))=(\nu ^2(dy)/ \nu ^1(dy)) \tilde \rho (y,x)$.}
\par  {\bf Proof.} (A). From the conditions imposed on $\mu ^{j,y}$
and $\nu ^j$ it follows that for each $\phi \in L(X\times Y, \pi ^j)$
due to Theorem 7.15\cite{roo} the following equality is accomplished
$$\int_{X\times Y}
\phi (x,y) \pi ^j(d(x,y))=\int_Y [ \int_X \phi (x,y) \mu ^{j,y}(dx)]
\nu ^j(dy),$$
 also $\rho (y,x)= \pi ^2(d(x,y))/ \pi ^1(d(x,y)) $
$\in L(X\times Y, \pi ^1)$, hence $\nu ^2(dy)/ \nu ^1(dy)$ $=[
\int_X \rho(y,x)\mu ^{1,y}(dx)] \in L(Y, \nu ^1)$. Further we modify
the proof of Theorem 1 \S 15\cite{sko}. Then $\tilde \rho (y,x)$
may be defined for $\nu ^2$-almost all $y$ by  $\tilde \rho (y,x)=
\rho (y,x)/ \int_X \rho(y,x) \mu ^{1,y}(dx) \in L(X, \mu ^{1,y})$.
\par  (B). Let $A \in Bco(X)\times Bco(Y)$, $A_y:=[y:$ $(x,y) \in A]$,
then $\pi ^j(A)=\int_Y \mu ^{j,y}(A_y) \nu ^j(dy)$. If $\| A\|
_{\pi ^1}=0$, then $\| A_y\| _{\mu ^{1,y}}N_{\nu ^1}(y)=0$ for each
$y \in Y$, consequently, $\| A\| _{\pi ^2}=0$, since $\nu ^2(dy)/
\nu ^1(dy) \in L(\nu ^1)$, $\mu ^{2,y}(dx)/ \mu ^{1,y}(dx) \in
L(\mu ^{1,y})$, $ \tilde \rho \in L(X\times Y, \pi ^1)$ and
Conditions $(i,ii)$ in \S 3.7 are satisfied. 
From this it follows that $\pi ^2(d(x,y))/\pi ^1(d(x,y))
\in L(X\times Y, \pi ^1)$, since $\nu ^2(dy)/ \nu ^1(dy) \in
L(X\times Y, \pi ^1)$ with $\sup_y \| X\| _{\mu ^{j,y}} < \infty $.
\par  {\bf 3.9. Definition.} For a Banach space $X$ 
over $\bf K$ an element $a \in X$
is called an admissible shift of a measure $\mu $ with values in
$\bf K_s$, if $\mu _a \ll \mu $, where $\mu _a(A)= \mu (S_{-a}A)$
for each $A$ in $E\supset Bco(X)$, 
$S_aA:=a+A$, $\rho (a,x):=\rho _{\mu }
(a,x):=\mu _a(dx)/ \mu (dx)$, $M_{\mu }:=[a \in X:$ $\mu _a \ll
\mu ]$ (see \S \S 2.1 and 2.31).
\par  {\bf 3.10. Properties of $M _{\mu }$ and $\rho $ from \S 3.9.}
\par {\bf I.} {\it The set $M _{\mu }$ is a semigroup by addition, $\rho (a+b,
x)=\rho (a,x)\rho (b,x-a)$ for each $a,b \in M_{\mu }$.}
\par  {\bf Proof.} For each continuous bounded $f: X\to \bf K_s$:
$\int_X f(x)\mu _{a+b}(dx)$ $=\int_X f(x+a+b)
\mu (dx)$ $=\int_X f(x+a)\rho (b,x) \mu (dx)$ $=\int_X f(x)
\rho (b,x-a) \rho (a,x) \mu (dx)$, since $\| X\| _{\mu } < \infty $
and $f(x) \rho (b, x-a) \in L(\mu )$, consequently, $\rho (b,x-a)
\rho (a,x)= \rho (a+b,x) \in L(\mu )$ as a function of $x$ and
$\mu_{a+b} \ll \mu $.
\par  {\bf II.} {\it If $a \in M_{\mu }$, $\rho (a,x) \ne 0$
$(mod$ $\mu )$, then $\mu _a \sim \mu $, $-a \in M_{\mu }$ and
$\rho (-a,x) =1/ \rho(a,x-a)$.}
\par  {\bf Proof.} For each continuous bounded
$f:X \to \bf K_s$: $\int_X f(x) \mu (dx)$ $=\int_X f(x) [\rho
(a,x)/ \rho (a,x)] \mu (dx)$ $=\int_X [\mu _a(dx)/ \mu (dx)]^{-1}
\mu _a(dx)$, since $\| X\| _{\mu} < \infty $, hence $\mu _a
\sim \mu $.
\par  {\bf III. } {\it If $\nu \ll \mu $ and $\nu (dx)/ \mu (dx)=
g(x)$, then $M_{\mu } \cap M_{\nu }=M_{\mu } \cap [a:$ $\mu ([x:$
$g(x)=0,$ $g(x-a)\rho _{\mu }(a,x) \ne 0])=0]$ and $\rho _{\nu }
(a,x)=[g(x-a)/g(x)]\rho _{\mu }(a,x)$ $(mod $ $\nu )$ for $a \in
M_{\mu }\cap M_{\nu }$.}
\par  {\bf Proof.} For each continuous bounded function $f: X\to
\bf K_s$: $a \in M_{\mu }$ and $\int_Xf(x+a) \nu (dx)$ $=\int_X
f(x)g(x-a)\rho _{\mu }(a,x) \mu (dx)$ such that $\mu ([
x:$ $ g(x)=0,$ $g(x-a) \rho _{\mu }(a,x) \ne 0])=0$ we have $\int_X
f(x+a) \nu (dx)$ $=\int_X f(x)[g(x-a)\rho _{\mu }(a,x)/
g(x)] \nu (dx)$, since $\| X\| _{\nu }+\| X\| _{\mu } < \infty $,
$N_{\nu }(x)=\inf_{Bco(X) \supset U \ni x} \sup_{y \in U}[$ $|g(y)|
N_{\mu }(y)]$, consequently, $a \in M_{\mu }\cap M_{\nu }$. If
$a \in M_{\mu } \cap M_{\nu }$, then 
$$\int_X f(x)\rho _{\nu }(a,x)
g(x)\mu (dx)=\int_X f(x)g(x-a)\rho _{\mu }(a,x) \mu (dx),$$
 consequently,
$\rho _{\nu }(a,x)g(x)$ $=g(x-a)\rho _{\mu }(a,x)$ $(mod$ $\mu )$
and $\mu ([x:$ $g(x)=0,$ $g(x-a)\rho _{\mu }(a,x) \ne 0])=0$.
\par  {\bf IV.} {\it If $\nu \sim \mu $, then $M_{\nu }=M_{\mu }$.}
\par {\bf V.} {\it  For $\mu $ with values in $\bf K_s$ and $X=K^m$,
$m \in \bf N$ a family $M_{\mu }$ with a distance function $r(a,b)=\|
\rho (a,x)-\rho (b,x)\| _{N_{\mu (x)}}$ is a complete pseudoultrametrizable
space.}
\par {\bf Proof. } Let $(a_n) \subset M_{\mu }$
be a Cauchy sequence relative to $r$, then $(a_n)$ is bounded in
$X$ by $\| *\| _X$, since for $\lim_{j \to \infty } \| a_{n_j}\| =
\infty $ and $r(a_{n_j}, a_{n_{j+1}})< p^{-j}$ for $f \in L(\mu )$
with a compact support we have $\| f(x+a_{n_j})-f(x+a_{n_1})\| _{N
_{\mu }} <1/p$. Then for $f$ with $\| f(x+a_{n_1})\| _{N_{\nu }}>1/2$
and $\| f\| _{N_{\nu }}=1$ we get a contradiction: $\lim_j \| f(x+a_{n_j})\|
_{N_{\mu }}>1/2-1/p \ge 0$. This is impossible because of compactness of
$supp(f)$. Therefore, $(a_n)$ is bounded, consequently,
there exists a subsequence $(a_{n_j})=:(b_j)$ weakly converging
in $X$ to $b\in X$. Therefore, $\theta _j(z)=\int_X \chi _e(z(x+b_j))
\mu (dx)$ $\chi _e(z(b_j)) \theta (z)=\int_X \chi _e(z(x)) \rho
(b_j,x) \mu (dx)$, $\lim _j z(b_j)=z(b)$ and $\lim _j \theta _j(z)=
\chi  _e(z(b)) \theta (z)$ for each $z \in X'$. From Theorem 9.20
\cite{roo} it follows that there is $\rho \in L(\mu )$ with $\lim_j \| \rho
(b_j,x)- \rho (x)\| _{N_{\mu }}=0$, since $L(\mu )$ is the
Banach space and $\mu _j$
corresponding to $\theta _j$ converges in the Banach space 
${\sf M}(X)$. Therefore, $\int_X \chi
_e(z(x)) \mu _b(dx)$ $=\int_X \chi _e(z(x)) \rho (x) \mu (dx)$
for each $z \in X'=K^m$, consequently, $\rho (x)=\mu _b(dx)/
\mu (dx)$.
\par   {\bf 3.11. Definition.} For a Banach space 
$X$ over $\bf K$ and a measure $\mu: Bco(X)\to \bf K_s$, 
$a \in X$, $\| a\| =1,$ a vector $a$
is called an admissible direction, if $a \in M_{\mu }^{\bf K}:=[z:$
$\| z\| _X=1,$ $\lambda z \in M_{\mu }$ and $\rho (\lambda z,x) \ne 0$
$(mod$ $\mu )$ (relative to $x$) and for each $\lambda \in {\bf K}]
\subset X$. Let $a \in
M_{\mu }^{\bf K}$ we denote by $L_1:={\bf K}a,$ $X_1=X\ominus L_1$, 
$\mu ^1$ and
$\tilde \mu ^1$ are the projections of $\mu $ onto $L_1$ and $X_1$
respectively, $\tilde
\mu =\mu ^1 \otimes \tilde \mu ^1$ be a measure on $Bco(X)$,
given by the the following equation
$\tilde \mu (A\times C)=\mu ^1(A)\tilde \mu ^1(C)$
on $Bco(L_1) \times Bco(X_1)$ and extended on $Bco(X)$, 
where $A\in Bco(L_1)$ and $C\in Bco(X_1)$.
\par {\bf 3.12. Definition and notes.} A measure $\mu : Bco(X)\to \bf K_s$ 
for a Banach space $X$
over $\bf K$ is called a quasi-invariant measure 
if $M_{\mu }$ contains a $\bf K$-linear
manifold $J_{\mu }$ dense in $X$. \
\par From \S 3.10 and Definition 3.11
it follows that $J_{\mu } \subset M_{\mu }^{\bf K}$. 
\par Let $(e_j:$
$j \in {\bf N})$ be orthonormal basis in $X$, $H=sp_{\bf K}(e_j:$ $j)$. 
We denote
$\Omega (Y)=[\mu | \mu $ is a measure with a finite total
variation on $Bco(X)$ and $H \subset J_{\mu }]$, where $Y=\bf K_s$. 
\par  {\bf 3.13. Theorem.} {\it If $\mu :Bf(Y) \to \bf F$ 
is a $\sigma $-finite measure on $Bco(Y)$, $Y$
is a complete separable ultrametrizable $\bf K$-linear subspace
such that $co(S)$ is nowhere dense in $Y$ for each compact
$S \subset Y$, where $\bf K$ and $\bf F$ are infinite non-discrete
non-Archimedean fields with multiplicative ultranorms
$|*|_{\bf K}$ and $|*|_{\bf F}$. Then from $J_{\mu }=Y$
it follows that $\mu =0$.}
\par  {\bf Proof.} Since $\mu $ is $\sigma $-finite, then
there are $(Y_j:$ $j \in H)\subset Bco(Y)$
such that $Y=\bigcup_{j\in H}
Y_j$ and $0< \| \mu |Bco(Y_j) \| \le 1$
for each $j$, where
$H \subset
\bf N$, $Y_j \cap Y_l= \emptyset $ for each $j \ne l$. If $card(
H)=\aleph _0$, then we define a function 
$f(x)=s^j/ \| Y_j\| _{\mu } $ for
$\mu $ with values in $\bf F$, where $s$ is fixed with
$0<|s|_F<1$, $s
\in \bf N$). Then we define a measure $\nu (A)=\int_A f(x) \mu (dx)$, $A
\in Bco(Y)$. Therefore, $\| Y\| _{\nu } \le 1$
and $J_{\nu }=Y$, since $f \in L(Y, \mu , {\bf F})$. 
Hence it is sufficient to consider
$\mu $ with $\| \mu \| \le 1$ and  $\mu (Y)=1$.
For each $n \in \bf N$ in view of the Radonian property of
$Y$ there exists a compact $X_n \subset Y$ such that $\| Y\setminus
X_n\| _{\mu }<s^{-n}$. In $Y$ there is a countable everywhere
dense subset $(x_j:$ $j \in {\bf N})$, hence $Y=\bigcup _{j
\in \bf N}B(Y,x_j,r_l)$
for each $r_l>0$, where $B(Y,x,r_l)=[y \in Y:$ $d(x,y) \le r_l]$,
$d$ is an ultrametric in $Y$, i.e. $d(x,z) \le max(d(x,y),$ $d(y,z))$,
$d(x,z)=d(z,x)$, $d(x,x)=0$, $d(x,y)>0$ for $x \ne y$ for each $x,y,z
\in Y$. Therefore, for each $r_l=1/l,$ $l\in \bf N$ there exists
$k(l)
\in \bf N$ such that $\| Y\setminus
X_{n,l}\| _{\mu }<s{-n-l}$ due to compactness of $Y_c=[y \in Y:$
$N_{\mu }(y) \ge c]$ for each $c>0$, where $X_{n,l}:=\bigcup_{j=1}
^{k(l)}B(Y,x_j,r_l)$, consequently, 
$\| Y\setminus X_n\| _{\mu } \le s^{-n}$ for $X_n:=
\bigcap_{l=1}^{\infty }X_{n,l}$. The subsets $X_n$ are compact,
since $X_n$ are closed in $Y$ and the metric $d$ on $X_n$ is
completely bounded and $Y$ is complete (see Theorems 3.1.2 and 4.3.29
\cite{eng}). Then $0< \| X\| _{\mu } \le 1$ for 
$\| Y\setminus X\| _{\mu }=0$ and for $X:=sp_K(\bigcup_{n=1}
^{\infty } X_n).$
\par  The sets $\tilde Y_n=co(Y_n)$ are nowhere dense in $Y$ for
$Y_n=\bigcup_{l=1}^nX_l$, consequently, $sp_{\bf K}Y_n$ are nowhere dense
in $Y$. Moreover, $(Y\setminus \bigcup_{n=1}^{\infty }Y_n) \ne \emptyset $
is dense in $Y$ due to the Baire category theorem (see 3.9.3 and 4.3.26
\cite{eng}).
Therefore, $y+X \subset Y\setminus X$ for $y \in Y\setminus X$ and from
$J_{\mu }
=Y$ it follows that $\| X\| _{\mu }=0$, since 
$\| y+X\| _{\mu }=0$ (see \S \S 2.32 and 3.12 above). 
Hence we get the contradiction, consequently, $\mu =0$.
\par  {\bf 3.14. Corollary.} {\it If $Y$ is a Banach space or a complete
countably-ultranormable infinite-dimensional over $\bf K$ space,
$\mu : Bco(Y) \to \bf K_s$, $\bf K$ and $\bf F$
are the same as in \S 3.13 and $J_{\mu }=Y$, then $\mu =0$.}
\par  {\bf Proof.} The space $Y$ is evidently complete and
ultrametrizable, since its topology is given by a countable
family of ultranorms. Moreover, $co(S)$ is nowhere dense in $Y$
for each compact $S$ in $Y$, since $co(S)=cl(S_{bc})$ is compact
in $Y$ and does not contain in itself any open subset in $Y$ due to
\S (5.7.5) in \cite{nari}.
\par {\bf 3.15. Theorem.} { \it Let $X$ be a separable Banach space over
a locally compact infinite field $\bf K$ with a nontrivial valuation
such that either ${\bf K}\supset \bf Q_p$ or $char ({\bf K})=p>0$. 
Then there are probability measures 
$\mu $ on $X$ with values in $\bf K_s$ $(s\ne p)$
such that $\mu $ are quasi-invariant 
relative to a dense $\bf K$-linear subspace
$J_{\mu }$}.
\par {\bf Proof.} Let $S(j,n):=p^{j}B({\bf K},0,1)
\setminus p^{j+1} B({\bf K},0,1)$
for $j\in \bf Z$ and $j\le n$,
$S(n,n):=p^nB({\bf K},0,1)$, $w'$ be the Haar measure on $\bf K$
considered as the additive group (see \cite{hew,roo}) with
values in $\bf K_s$ for $s\ne p$.
Then for each $c>0$ and $n \in {\bf N}$ there are measures
$m$ on $Bf({\bf K})$ such that
$m(dx)=f(x)v(dx)$, $\mid f(x)\mid >0$ for each $x \in \bf K$
and $\mid m(p^{n}B({\bf K},0,1))-1\mid <c$,
$m({\bf K})=1$, $\mid m\mid (E)\le 1$ for each $E\in Bco({\bf K})$ ,
where $v=w'$, $v(B({\bf K},0,1))=1$. Moreover,
we can choose $f$ such that a density $m_{a}(dx)/m(dx)=:d(m;a,x)$
be continuous by
$(a,x)\in \bf K^{2}$  and for each $c'>0$,  $x$ and
$\mid a\mid \le p^{-n}$ : $\mid d(m;a,x)-1\mid <c'.$
Let $f|_{S(j,n)}:=a(j,n)$ be locally constant,
for example, $a(j,n)=(1-s)(1-1/p)s^{2n-1-j}p^{-n}$ for $j<n$, $a(n,n)=
(1-s^{-n})p^{-n}$. Then taking $f+h$ and using $h(x)$ 
with $0<\sup_x|h(x)/f(x)|=c" \le 1/s^n$ we get
$\mid y_{a}(dx)/y(dx)\mid =\mid m_{a}(dx)/m(dx)\mid ,$
where $y(dx)=(f+h)(x)m(dx)$.
\par  Let $\{ m(j;dx) \}$ be a family of measures on
$\bf K$ with the corresponding sequence $\{ k(j) \}$ such that
$k(j)\le k(j+1)$ for each $j$ and
$\lim_{i\to \infty }k(i)= \infty $, where
$m(j;dx)$ corresponds to the partition $[S(i,k(j))]$. The Banach space 
$X$ is isomorphic with $c_0
(\omega _0,{\bf K})$ \cite{roo}. It has the orthonormal basis
$\{ e_j:j=1,2,... \} $ and the projectors
$P_jx=(x(1),...,x(j))$ onto $\bf K^j$, where $x=x(1)e_1+x(2)e_2+...$.
Then there exists a cylindrical measure $\mu $ generated by a consistent family
of measures $y(j,B)=b(j,E)$ for $B=P^{-1}_jE$ and $E\in Bf({\bf K}^j)$ 
\cite{boui,dal}
where $b(j,dz)=\otimes [m(j;dz(i)): i=1,...,j], z=(z(1),...,z(j))$.
Let $L:=L(t,t(1),...,t(l);l):=\{ x: x \in X \mbox{ and } \mid x(i)\mid
\le p^{a}, a=-t-t(i) \mbox{ for } i=1,...,l, \mbox{ and } a=-k(j)
\mbox{ for } j>l \} $, then $L$ is compact in $X$, since $X$ is
Lindel\"of and $L$ is sequentially compact \cite{eng}.
Therefore, for each $c>0$ there exists
$L$ such that $\| X\setminus L\|
_{\mu }<c$ due to
the choice of $a(j,n)$.
\par In view of the Prohorov theorem  for measures with values in
$\bf K_s$ 7.6(ii)\cite{roo} and due to Lemma
2.3 $\mu $ has the countably-additive extension on $Bf(X)$, 
consequently, also on the complete $\sigma $-field
$Af(X,\mu )$ and $\mu $ is the Radon measure.
\par Let $z' \in sp_{\bf K} \{ e_j:j=1,2,... \} $
and $z"= \{ z(j): z(j)=0
\mbox{ for } j\le l \mbox{ and } z(j)\in S(n,n), j=1,2,...,
n=k(j) \} $, $l \in {\bf N}$, $z=z'+z"$.
Now take the restriction of $\mu $ on $Bco(X)$.
In view of Theorems 2.19, 3.5 above and also I.1.4, II.4.1\cite{roo}
there are $m(j;dz(j)$ such that $\rho _{\mu }(z,x)=\prod \{ d(j;z(j),
x(j)): j=1,2,...\} =
\mu _z(dx)/\mu (dx)\in L((X,\mu ,Bco(X)),{\bf K_s})$
for each such $z$ and $x\in X$, where $d(j;*,*)=
d(m(j;*),*,*)$ and $\mu _z(X)=\mu (X)=1$.
\par {\bf 3.16. Note.} For a given $m=w'$ (see above)
new suitable
measures may be constructed, if to use images of measures
$m^{g}(E)=m(g^{-1}(E))$ such that for a diffeomorphism
$g\in Diff^1({\bf K})$ (see \S A.3) we have $m^{{g}^{-1}}(dx)/m(dx)=
\mid (g'(g^{-1}(x))\mid _{\bf K}$, where $|*| _{\bf K}= mod_{\bf K}(*)$ 
is the modular function of the field $\bf K$ associated with
the  Haar meassure on $\bf K$, at the same time  
$|*| _{\bf K}$ is the multiplicative norm in $\bf K$
consistent with its uniformity \cite{wei}.
Indeed, for $\bf K$ and $X=\bf K^j$ with $j \in \bf N$ and 
the Haar measure $v=w'$ on $X$, $v_X:=v$
with values in $\bf K_s$ for $s\ne p$
and for a function $f \in L(X,v,{\bf K_s})$
we have: $\int_{g(A)} f(x)v(dx)$ $=\int_Af(g(y)) |g'(y)|_{\bf K} v(dy)$, where
$mod_{\bf K}(\lambda )v(dx):=v(\lambda dx)$, $\lambda \in \bf K$, 
since $v(B(X,
0,p^n))\in \bf Q$, $N_v(x)=1$ for each $x \in X$, consequently,
from $f_k \to f$ in $L(g(A),v,{\bf K_s})$ whilst $k \to \infty $ it
follows that
$f_k(g(x)) \to f(g(x))$ in $L(A,v,{\bf K_s})$, where
$f_k$ are locally constant, $A$ is compact and open in $X$.
\par Henceforward, quasi-invariant measure $\mu $ on
$Bco(c_0(\omega _0,{\bf K}))$
constructed with the help of projective limits or sequences
of weak distributions 
of probability measures $(\mu _{H(n)}:$ $n)$ are considered,
for example, as in Theorem 3.15 such that
\par  (i) $\mu _{H(n)}(dx)=f_{H(n)}(x)v_{H(n)}(dx)$, $dim_{\bf K}H(n)
=m(n)< \aleph _0$ for each $n \in \bf N$, where 
$f_{H(n)} \in L(H(n),v_{H(n)},{\bf K_s})$,
$H(n) \subset H(n+1) \subset ...$, $cl(\bigcup_nH(n)=c_0(
\omega _0,{\bf K})$, if it is not specified in another manner.
\par For probability quasi-invaraitn measure
with values in $\bf K_s$,
if shifts $x \mapsto x+y$ by $y \in H(n)$ are continuous from $H(n)$ to
${\sf M}(H(n))$ (see \S 2.1), that is, $y \to \mu ^y_{H(n)}$, where
$\mu _{H(n)}(y+A)=:
\mu ^y_{H(n)}(A)$ for $A \in Bco(H(n))$, then due to Theorem 8.9
\cite{roo} $\mu _{H(n)}$ satisfies (i).
\par  As will be seen below such measures $\mu $ are quasi-invariant
relative to families of the cardinality ${\sf c}=card ({\bf R})$
of linear and non-linear transformations $U: X\to X$. Moreover, for each
$V$ open in $X$  we have $\| V\| _{\mu }>0$,
when $f_{H(n)}(x) \ne 0$ for each $n \in \bf N$ and $x \in H(n)$.
\par  Let $\mu $ be a probability 
quasi-invariant measure satisfying (i) and  $(e_j:$ $j)$
be orthonormal basis in $ M_{\mu }$, 
$H(n):=sp_{\bf K}(e_1,...,e_n)$, we denote by
\par  $\hat \rho
_{\mu }(a,x)=\hat \rho (a,x)=\lim_{n \to \infty } \rho ^n(P_na,P_nx)$,
\par $\rho ^n(P_na,P_nx):=f_{H(n)}(P_n(x-a))/f_{H(n)}(P_nx)$ for each
$a$ and $x$ for which this limit exists and $\hat \rho (a,x)=0$
in the contrary case, where $P_n: X\to H(n)$ are chosen 
consistent projectors. Let
$\rho (a,x)=\hat \rho (a,x)$, if $\mu _a(X)=\mu (X)$ and
$\hat \rho (a,x)\in L(X, \mu ,{\bf K_s})$ as a function by $x$ and
$\| X\| _{N_{\nu }}=1$,
where $\nu (dx):=\hat \rho (a,x) \mu (dx)$, $\rho (a,x)$ is not
defined when $\mu _a(X)=\mu (X)$ or $\| X\| _{N_{\nu }} \ne 1$,
this condition of the equality to $1$ may be satisfied, for example,
for continuous $f_{H(n)}$ with continuous
$\hat \rho (a,x)\in L(\mu )$ by $x$ for each given $a$, if
$\lim_n\rho ^n(a,x)$ converges uniformly by $x$.
If for some another basis $(\tilde e_j:$ $j)$ and $\tilde \rho $
is accomplished
\par  (ii) $\| X\setminus S\| _{\mu }=0$,
then $\rho (a,x) $ is called regularly dependent from a basis, where
$S:=\bigcap_{a \in M_{\mu }}[x:$ $\rho
(a,x)=\tilde \rho (a,x)])$.
\par  {\bf 3.17. Lemma.} {\it Let $\mu $ be a probability measure,
$\mu : Bco(X)\to \bf K_s$, $X$ be a Banach space over $\bf K$,
suppose that for each basis $(\tilde e_j:$ $j)$ in
$M_{\mu }$ a quasi-invariance factor $\tilde \rho $ satisfies the following
conditions: 
\par $(1)$ if $\tilde \rho (a_j,x)$, $j=1,...,N$, are defined for
a given $x \in X$ and for each $\lambda _j \in \bf K$ then a function
$\tilde \rho (\sum_{j=1}^n \lambda _ja_j,x)$ is continuous by
$\lambda _j$, $j=1,...,N$; 
\par $(2)$ there exists an increasing sequence of subspaces
$H(n) \subset M_{\mu }$, $cl(\bigcup _nH(n))=X$, with projectors
$P_n: X \to H(n)$, $B \in Bf(X)$, $\| B\| _{\mu }=0$
such that $\lim_{n\to \infty }\tilde \rho (P_na,x)=\tilde \rho (a,x)$
for each $a \in M_{\mu }$ and $x \notin B$ for which is defined
$\rho (a,x)$. Then $\rho (a,x)$ depends regularly from the basis.}
\par  {\bf Proof.} There exists a subset $S$ dense in each
$H(n)$, hence $\| B'\| _{\mu }=0$ for
$B'=\bigcup_{a \in S}[x:$ $\rho (a,x)\ne \tilde \rho (a,x)$].
From (1) it follows that $\tilde \rho (a,x)=\rho (a,x)$ on each $H(n)$
for $x \notin B'$. From $sp_{\bf K}S \supset H(n)$ and (2) it follows that
$\rho (a,x)=\tilde \rho (a,x)$ for each $a \in M_{\mu }$ and $x \in
X\setminus (B'\cup B)$, consequently, Condition 3.16.(ii) is satisfied, since
from $\rho (a,x) \in L(\mu )$ it follows that $\tilde \rho (a,x)
\in L(\mu )$ as the function by $x$.
\par {\bf 3.18. Lemma.} \
{\it If a probability quasi-invariant measure $\mu : Bco(X)\to \bf K_s$
satisfies Condition 3.16.(i), then there exists a compact operator
$T: X\to X$ such that $M_{\mu } \subset (TX)^{\sim }$, where $X$ is
the Banach space over $\bf K$.}
\par {\bf Proof.} Products of dense measures are dense measures
due to Theorem 7.28\cite{roo}, whence for $\mu _{H(n)}(dx)=
\bigotimes _{j=1}^{m(n)}\mu _{{\bf K}e(j)}(dx_j)$ is accomplished
$N_{\mu _{H(n)}}
(x)=\prod_{j=1}^{m(n)} N_{\mu _{{\bf K}e(j)}}(x_j)$, where
$x=(x_1,..,x_{m(n)})$,
$x_j \in \bf K$. From Theorem 7.6\cite{roo} and Lemma I.2.5 it follows that
for each $1>c>0$ there are $R_j=R_j(c)$ with $[x_j:$
$N_{\mu _{{\bf K}e(j)}}(x_j) \ge
c] \subset B({\bf K},0,R_j)$ and $\lim_{j \to \infty }R_j=0$. Choosing
$c=c(n)=s^{-n}$, $n \in \bf N$ and using $\prod_{j=1}^{\infty }
=\epsilon _j=0$ whilst $0<\epsilon _j<c<1$ for each $j$ we get
that there exists a sequence $[r_j:$ $j]$ for which $card[
j:$ $|a_j|>r_j]<\aleph _0$ for each $a \in M_{\mu }$, since
$[x\in X:$ $|x_j|\le r_j$ for all $j]$ is a compact subgroup in
$X$, where $a=(a_j:$ $j)$, $a_j \in \bf K$,
$r_j>0$, $\lim_j r_j=0$. Therefore, $M_{\mu } \subset (TX)^{\sim }$
for $T=diag(T_j:$ $j)$ and $|T_j| \ge r_j$ for $j \in \bf N$.
\par  {\bf 3.19.} Let $X$ be a Banach space over $\bf K$,
$|*|_{\bf K}=mod_{\bf K}(*)$, $U: X\to X$ be an invertible linear operator,
$\mu : Bco(X)\to 
\bf K_s$ be a probability quasi-invariant measure.
\par The uniform convergence of a (transfinite) sequence of functions on
$Af(V,\nu )$-compact subsets of a topological space
$V$ is called the Egorov condition, where $\nu $ is a measure on  $V$.
\par {\bf Theorem.} {\it Let pairs $(x-Ux,x)$
and $(x-U^{-1},x)$ be in $dom(\tilde \rho (a,x))$, where $dom(f)$ denotes
a domain of a function $f$, $\tilde \rho (x-Ux,x)\ne 0$, $\tilde
\rho (x-U^{-1}x,x)\ne 0$ $(mod $ $\mu )$ and $\mu $ satisfies 
Condition 3.16.(i), also $\tilde \rho
(\tilde P_n(x-Ux),x)=:\tilde \rho _n(x) \in L(\mu )$ and $\tilde \rho (
\tilde P_n(x-U^{-1}x,x)=:\bar \rho _n(x) \in L(\mu )$ converge
uniformly on $Af(X, \mu )$-compact subsets in $X$ such that there
exists $g \in L(\mu )$ with $|\tilde \rho _n(x)| \le |g(x)|$ and
$|\bar \rho _n(x)| \le |g(x)|$ for each $x \in X$ and each
projectors $\tilde P_nX\to \tilde H(n)$ with $cl(\bigcup_n
\tilde H(n))=X$, $\tilde H(n) \subset \tilde H(n+1)\subset ...$,
that is, Egorov conditions are satisfied for $\tilde \rho _n$ and
$\bar \rho _n$.
Then $\nu \sim \mu $ and
$$(i) \mbox{ }\nu (dx)/ \mu (dx)=|det(U)|_{\bf K}\tilde \rho (x-U^{-1}x,x),$$
if $\rho $ depends regularly from the base, then $\tilde \rho $ may be
substituted by $\rho $ in formula $(i)$, where $\nu (A):=
\mu (U^{-1}A)$ for each $A\in Bco(X)$.}
\par  {\bf Proof.} In view of Lemma 3.18 there
exists a compact operator $T: X\to X$ such that $M_{\mu }
\subset (TX)^{\sim }$, consequently, $(U-I)$ is a compact operator,
where $I$ is the identity operator. From the invertibility of
$U$ it follows that $(U^{-1}-I)$ is also compact, moreover, there
exists $det(U) \in \bf K$. Let
$g$ be a continuous bounded function, $g: \tilde H(n) \to
\bf K_s$, whence $\int_X\phi (x)\nu (dx)$ $=\int_{\tilde H(n)}
g(x)[f_{H(n)}(U^{-1}x)/f_{\tilde H(n)}(x)]|det(U_n)|_{\bf K} \mu _{\tilde
H(n)}(dx)$, for $\phi (x)=g(\tilde P_nx)$, where subspaces exist such
that $\tilde H(n) \subset X$, $(U^{-1}-I)\tilde H(n) \subset \tilde H(n)$,
$cl(\bigcup_n \tilde H(n))=X$, $U_n:=\hat r_n(U)$, $r_n=\tilde P_n:
X \to \tilde H(n)$ (see \S \S I.3.8 and II.3.16), $\tilde H(n) \subset \tilde
H(n+1)\subset ...$ due to compactness of $(U-I)$. In view of
the Lebesgue theorem due to fulfillment of the Egorov conditions
for $\tilde \rho _n$ and $\bar \rho _n$, see \S 7.6 \cite{mon} or
\S 7.F \cite{roo} $J_m=J_{m,\rho }$, since
$\tilde \rho (x-U^{-1}x,x)\in L(\mu )$, where $J_m:=\int_Xg(\tilde P_mx)
\nu (dx)$ and $J_{m,\rho }:=\int_Xg(\tilde P_mx)\tilde \rho (x-U^{-1}x,x)
|det(U)|_{\bf K} \mu (dx)$. Indeed, there exists $n_0$ such that
$|u(i,j)-\delta _{i,j}| \le 1/p$ for each $i$ and $j>n_0$,
consequently, $|det(U_n)|_{\bf K}=|det(U)|_{\bf K}$ for each
$n>n_0$. Then due to Condition 3.16.(i) 
and the Egorov conditions (see also \S 3.3) 
there exists $\lim_{n\to \infty }[\mu _{\tilde
H(n)}(d\tilde P_nx)/\nu _{\tilde H(n)}(d\tilde P_nx)]$
$=\mu (dx)/ \nu (dx)$ $(mod $ $\nu )$. Further analogously to the proofs
of Theorems 1 and 2 \S 25\cite{sko}.
\par {\bf 3.20.} Let $X$ be a Banach space over $\bf K$,
$|*|_{\bf K}=mod_{\bf K}(*)$ with a probability quasi-invariant 
measure $\mu : Bco(X)\to \bf K_s$ and Condition 3.16.(i) be satisfied,
also let $U$ fulfils the following conditions:
$$(i)\mbox{ }U(x) \mbox{ and } U^{-1}(x)
\in C^1(X,X) (see \S I.A.3);$$
$$ (ii)\mbox{ } (U'(x)-I)\mbox{ is compact for each }x\in X;$$
$$(iii)\mbox{ }(x-U^{-1}(x)) \mbox{ and }(x-U(x)) \in J_{\mu }
\mbox{ for }\mu -\mbox{a.e. }x\in X;$$
$$ (iv)\mbox{ for }\mu \mbox{-a.e. }x \mbox{ pairs }
(x-U(x);x) \mbox{ and }(x-U^{-1}(x);x)$$
are contained in a domain of $\rho (z,x)$ such that
$\rho (x-U^{-1}(x),x)\ne 0$, $\rho (x-U(x),x)\ne 0$ $(mod$ $\mu )$;
$$ (v)\mbox{ }\| X\setminus S'\| _{\mu }=0,$$
where $S':=([x: \mbox{ }\rho (z,x) \mbox{ is defined and continuous by }
z\in L ])$
for each finite-dimensional $L\subset J_{\mu }$;
$$ (vi)\mbox{ there exists }S \mbox{ with }
\| S\| _{\mu }=0\mbox{ and for each }$$ $x\in X\setminus S$
and for each $z$ for which there exists $\rho (z,x)$ satisfying the
following condition:
$\lim_{n \to \infty }\rho (P_nz,x)=\rho (z,x)$ and the convergence
is uniform for each finite-dimensional $L\subset J_{\mu }$ by $z$ in
$L\cap [x\in J_{\mu }:$ $\mid x\mid \le c]$, where $c>0$,
$P_n: X\to H(n)$  are projectors onto
finite-dimensional subspaces $H(n)$ over $\bf K$ such that
$H(n)\subset H(n+1)$ for each $n \in \bf N$ and $cl \cup \{ H(n): n \}=X$;
$$ (vii)\mbox{ there exists }n \mbox{ for which for all }j>n
\mbox{ and }
x\in X  \mbox{ mappings}$$
$V(j,x):=x+P_j(U^{-1}(x)-x)$ and $U(j,x):=
x+P_j(U(x)-x)$ are invertible and
$\lim_j \mid \det U'(j,x)\mid =\mid \det U'(x)\mid $,
$\lim_j \mid \det V'(j,x)
\mid =1/ \mid \det U'(x)\mid $ with the Egorov convergence
in (vi) by $z$
for $\rho (P_nz,x)$ and in (vii) by $x$ for $|det (U'(j,x))|$ and
$|det(V'(j,x))|$ for $\mu $ with values in $\bf K_s$.
\par {\bf Theorem.} {\it The measure 
$\nu (A):=\mu (U^{-1}(A))$ is equivalent to $\mu $ and
$$ (i)\mbox{ }\nu (dx)/ \mu (dx)=\mid \det U'(U^{-1}(x))\mid _{\bf K}
\rho (x-U^{-1}(x),x).$$}
\par {\bf Proof.} The beginning of the proof is analogous to that of
\S I.3.25.
Due to Conditions (vi, vii) we get $\lim_n \rho (x-V_n^{-1}x,x)=
\rho (x-U_1^{-1}x,x)$ in $L(\mu )$
due to the Egorov conditions. Then $J_1=J_{1,\rho }$ due to the Lebesgue
convergence theorem,
where $J_1=\int_Xf(U_1x) \mu (dx)$, $J_{1,\rho }:=\int_Xf(x)\rho (x-
U_1^{-1}x,x)|det U_1|_{\bf K}\mu (dx)$ for continuous bounded function
$f:X\to \bf K_s$. Analogously for $U_1^{-1}$
instead of $U_1$. Using instead of $f$ the function $\bar \Phi (U_1^{-1}
x):=f(x)\rho _{\mu }(x-U_1^{-1}x,x)$ and Properties 3.10 we get
that $\rho _{\mu }(U_1x-x,U_1x)\rho _{\mu }(x-U_1x,x)=1$ $(mod$ $\mu )$.
Therefore, for $U=U_1U_2$ with diagonal $U_1$ and upper triangular $U_2$
and lower triangular $U_3$ operators with finite-dimensional over
$\bf K$ subspaces $(U_j-I)X$, $j=1,2,3$, the following
equation is accomplished
$\int_Xf(Ux)\mu (dx)$ $=\int_X f(x)\rho _{\mu }(x-U^{-1}x,x)
|det U|_{\bf K} \mu (dx)$. If $(S^{-1}U-I)X=L$ [or $(U^{-1}S-I)X=L$],
then from the
decomposition given in (I) $U=SU_2U_1U_3,$ we have $(U_j-I)X=L$,
[or $(U_j^{-1}-I)X=L$
respectively], $j=1,2,3$ due to formulas from \S I.A.1, since
corresponding non-major minors are equal to zero.
\par If $U$ is an arbitrary linear operator satisfying
the conditions of this theorem, then from (iv-vi) and (I, II)
for each continuous bounded function $f: X \to \bf K_s$
we have $J=J_{\rho }$, where $J:=\int_Xf(U(x))\mu (dx)$ and $J_{\rho }:=
\int_Xf(x)\rho _{\mu }(x-U^{-1}(x),x)|det U|_{\bf K}\mu (dx)$. 
Analogously for
$U^{-1}$, moreover, $\rho (x-U^{-1}(x),x)|det U|_{\bf K}=:h(x)\in L(\mu )$,
 $h(x)\ne 0$ $(mod$ $\mu )$, since there exists
$det U$.
\par Suppose $U$ is polygonal (see \S I.3.25).
Then $U^{-1}$ is also polygonal, $U'(x)=V(j)$ for $x \in Y(j)$ and
$\int_X f(a(i)+V(i)x)\mu (dx)=\int f(a(i)+x)\rho _{\mu }(x-V^{-1}(i)x,x)$
$\times |det(V(i))|_{\bf K} \mu (dx)$
for each continuous bounded
$f: X \to \bf K_s$ and each $i$. From $a(j) \in M_{\mu }$
and \S 3.10 we get $\int_Xf(a(j)+V(j)x) \mu (dx)$ 
$=\int_X f(x)\rho (x-V(j)^{-1}
(x-a(j)),x) |det V(j)|_{\bf K} \mu (dx)$. Let $H_{k,j}:=[x\in X:$ $V(k)^{-1}x
=V(j)^{-1}x]$, assume without loss of generality that $V(k) \ne V(j)$ or
$a(k) \ne a(j)$ for each
$k \ne j$, since $Y(k) \ne Y(j)$ (otherwise they may be united).
Therefore, $H_{k,j} \ne X$. If $\| H_{k,j}
\| _{\mu }>0$, then from $X \ominus H_{k,j} \supset \bf K$ it follows
that $M_{\mu }
\subset H_{k,j}$, but $cl(H_{k,j})=H_{k,j}$ and $cl(M_{\mu })=X$. This
contradiction means that $\| A\| _{\mu }=0$,
where $A=[x:$ $V(k)^{-1}(x-a(k))=V(j)^{-1}(x-a(j))]$. Then
$\int_Xf(U(x)) \mu (dx)$ $=\int_Xf(x)\rho (x-U^{-1}(x),x)$
$|det U'(x)|_{\bf K}^{-1}\mu (dx)$.
\par  Then as in \S I.3.25.(V) for the construction
of the sequence $\{ U(j,*) : j \}$ it is sufficient to construct a sequence
of polygonal functions $\{ a(i,j;x)\} $, that is $a(i,j;x)=l_k(i,j)(x)+a_k$
for $x \in Y(k)$, where $l_k(i,j)$ are linear functionals,  $a_k\in \bf K$,
$Y(k)$ are closed in $X$, $Int(Y(j))\cap Int(Y(k))=\emptyset $ for each
$k \ne j$, $\bigcup_{k=1}^m Y(k)=X$, $m<\aleph _0$. For each $c>0$
there exists $V_c \subset X$ with  $\| X\setminus V_c\| <c$,
 the functions $s(i,j;x)$ and
$(\bar \Phi ^1 s(i,j; *))(x,e(k),t)$
are equiuniformly continuous (by $x \in V_c$ and by $i, j, k \in {\bf N}$)
on $V_c$. Choosing $c=c(n)=s^{-n}$
and using $\delta $-nets in $V_c$ we get a sequence of polygonal
mappings $(W_n:$ $n)$ converging by its matrix elements by Egorov
in the Banach space $L(X,\mu ,{\bf K_s})$,
from Condition (i) it follows that it may be chosen equicontinuous for
matrix elements $s(i,j;x)$, $ds(i,j;x)$
and $s(i,P_jx)$ by $i, j$ (the same is true for $U^{-1}$).
\par  Then calculating integrals as above for $W_n$ with functions
$f$, using the Lebesgue convergence theorem we get
the equalities analogous to written in
\S I.3.25.(III) for $J$ and $J_{\rho }$ of the general form. From
$\nu (dx)/\mu (dx)\ne 0$
$(mod$ $\mu )$ and \S 2.19 we get the statement of the theorem.
\par   {\bf 3.21. Examples.} Let $X$ be a Banach space 
over the field $\bf K$ with the valuation group 
$\Gamma _{\bf K}=\Gamma _{\bf Q_p}$.
We consider a diagonal compact operator
$T=diag(t_j:$ $j\in {\bf N})$ in a fixed orthonormal basis
$(e_j: j )$ in $X$ such that $ker T:=T^{-1}0=\{ 0 \}$. Let
$\nu '_j(dx_j)=C'(\xi _j)s^{-q \min (0,ord_p((x_j-x^0_j)/\xi _j))}
v(dx_j)$ for the Haar measure $v: Bco({\bf K})\to \bf Q_s$, 
then $\nu '_j(Bf({\bf K})) \subset
\bf C_s$. We choose constant functions $C'(\xi _j)$
such that $\nu '_j$ be a probability measure, where
$x^0=(x^0_j:$ $j) \in X$, $x=(x_j:$ $j) \in X$, $x_j \in \bf K$.
\par With the help of products
$\bigotimes_j \nu '_j(dx_j)$ as in \S 3.15 we can construct a probability
quasi-invariant measure $\mu ^T$ on $X$ with values in $\bf C_s$, 
since $cl(TX)$ is compact in $X$ and $sp_{\bf K}(e_j:$ $j)=:H
\subset J_{\mu }$. From $\bigcap
_{\lambda \in B({\bf K},0,1)\setminus 0}cl(\lambda TX)=\{ 0\}$ we may infer
that for each $c>0$ there exists a compact $V_c(\lambda )
\subset X$ such that $\| X
\setminus V_c(\lambda )\| _{\mu }<c$ and $\bigcap_{\lambda \ne 0}
V_c(\lambda )=\{ 0\}$, consequently, $\lim_{|\lambda | \to 0}
\int_Xf(x) \mu ^{\lambda T}(dx)=f(0)=\delta _0(f)$, hence
$\mu ^{\lambda T}$
is weakly converging to $\delta _0$ whilst $|\lambda | \to 0$ for the
space of bounded continuous functions $f: X \to \bf C_s$.
\par From Theorem
3.6 we conclude that from $\sum_{j=1}^{\infty }|y_j/\xi _j|^q_p < \infty $
it follows $y \in J_{\mu ^T}$. Then for a linear transformation
$U: X\to X$ from $\sum_j|\tilde e_j(x-U(x))/\xi _j|_p^q < \infty $
it follows that $x-U(x) \in J_{\mu }$ and a pair $(x-U(x),x) \in dom(\rho
(a,z))$. Moreover, for $\rho $ corresponding to $\mu ^T$ conditions
(v) and (vi) in \S 3.20 are satisfied. Therefore, for such $y$ and
$S \in Af(X,\mu )$ 
a quantity $|\mu (ty+S)- \mu (S)|$ is of order of smallness $|t|^q$
whilst $t \to 0$, hence they are pseudo-differentiable of order $b$
for $0<Re(b)<q$ (see also \S 4 below).
\par It is interesting also to discuss a way of solution of one 
problem formulated in \cite{khen} that there does not exist
a $\sigma $-additive $\bf Q_p$-valued measure with values in  
$X$ over $\bf Q_p$ such that it would be an analog
of the classical Gaussian measure. 
In the clasical case this means in particular
a quasi-invariance of a measure relative to shifts
on vectors from a dense subspace.
We will show, that on a Banach space $X$ over ${\bf K}\supset \bf Q_p$ 
for each prime number $p$ there is not a $\sigma $-additive $\mu \ne 0$
with values in $\bf K_p$ such that it is quasi-invariant relative to 
shifts from a dense subspace.
Details can be lightly extracted from the results given above.
Let on $(X, Bco(X))$ there exists such $\mu $. 
With the help of suitable compact operators a cylindrical measure
on an algebra of cylindrical subsets of $X$ generates
quasi-invariant measures, so we can suppose that $\mu $ is quasi-invariant.
Then it produces a sequence of a finite-dimensional distribution
$\{ \mu _{L_n}: n\in {\bf N} \}$ analogously to \S 2 and \S 3,
where $L_n$ are subspaces of $X$ with dimensions over $\bf K$ 
equal to $n$.
Each measure $\mu _{L_n}$ is $\sigma $-additive.
From the quasi-invariance of $\mu $ it follows, that $L_n$ 
can be chosen such that $\mu _{L_n}$ are quasi-invariant
relative to the entire $L_n$. But in view of 
Chapters 7-9 \cite{roo} and \cite{sch1} 
for measures with values in $\bf K_p$
(see also Proposition 11 from \S VII.1.9 \cite{boui}) 
this means that $\mu _{L_n}$ is equivalent to the Haar measure
on $L_n$ with values in $\bf K_p$. The space $L_n$ as the additive group
can be considered over $\bf Q_p$,
moreover, for each continuous linear functional $\phi : 
{\bf K_p}\to \bf Q_p$ considered as the finite-dimensional Banach space
over $\bf Q_p$ the measure $\phi \circ \mu _{L_n}(*)$ 
is non-trivial for some $\phi $.
Consequently, on $L_n$ there would be the Haar measure 
with values in $\bf Q_p$, but this is impossible due to Chapter 9 
in \cite{roo}, since $L_n$ is not the $p$-free group. 
We get the contradiction, that is, such $\mu $ does not exist.
\section{ Pseudo-differentiable measures.}
\par {\bf 4.1. Definition.} A function $f: {\bf K}\to 
\bf \Lambda _s$ is called
pseudo-differentiable of order $b$, if there exists
the following integral:
$PD(b,f(x)):=\int_{\bf K}[(f(x)-f(y)) \times g(x,y,b)]
dv(y)$. We introduce the following notation $PD_c(b,f(x))$ for
such integral by $B({\bf K},0,1)$ instead of the entire $\bf K$.
Where $g(x,y,b):=s^{(-1-b) \times
ord_p(x-y)}$ with the corresponding Haar measure $v$
with values in $\bf K_s$, where $b \in {\bf C_s}$ and 
$|x|_{\bf K}=p^{-ord_p(x)}$,
$\bf C_p$ denotes the field of complex numbers with the 
non-Archimedean valuation extending that of $\bf Q_p$, 
$p^{-ord_p(\zeta )}:=|\zeta |_{\bf K},$  
$\bf \Lambda _p$ is a spherically complete
field with a valuation group $ \{ |x|:$ $0\ne x \in {\bf \Lambda _p} \} =
(0,\infty )\subset \bf R$ such that ${\bf C_p}\subset \bf \Lambda _p$
\cite{diar,roo,sch1,wei}. 
\par A quasi-invariant measure $\mu $ on $X$
is called pseudo-differentiable for $b \in \bf C_s$,
if there exists $PD(b,g(x))$ for $g(x):=\mu (-xz+S)$ for each
$S \in Bco(X)$ $\| S \| _{\mu }< \infty $
and each $z \in J^b_{\mu }$, where $J^b_{\mu }$ is a $\bf K$-linear
subspace dense in $X$. For a fixed $z \in X$ such measure is called
pseudo-differentiable along $z$.
\par  For a one-parameter subfamily of 
operators $B({\bf K},0,1)\ni t \mapsto
U_t: X \to X$ quasi-invariant measure $\mu $ is called 
pseudo-differentiable for $b \in \bf C_s$,
if for each $S$ the same as above there exists $PD_c(b,g(t))$ for a function
$g(t):=\mu (U_t^{-1}(S)$, where $X$ may be also a topological group
$G$ with a measure quasi-invariant relative to a dense subgroup
$G'$ (see \cite{lu3,luum983,lubp2}).
\par {\bf 4.2.} Let $\mu $, $X$, and $\rho $ be the same as in Theorem 3.15
and $F$ be a non-Archimedean Fourier transform defined in \cite{vla,roo}.
\par {\bf Theorems.} {\it $(1) $ $g(t):=\rho (z+tw,x)j(t) 
\in L({\bf K},v,{\bf K_s})=:V$ for $\mu $ and the Haar measure
$v$ with values in $\bf K_s$, where
$z$ and $w \in J_{\mu }$, $t \in \bf K$, 
$j(t)$ is the characteristic function
of a compact subset $W \subset \bf K$. In general, may be
$k(t):=\rho (z+tw,x) \notin V$.
\par  $(2)$ Let $g(t)=\rho (z+tw,x)j(t)$
with clopen subsets $W$ in $\bf K$. Then there are $\mu $,
for which there exists
$PD(b,g(t))$ for each $b \in {\bf C_s}$. If
$g(t)=\rho (z+tw,x)$, then there are probability measures $\mu $,
for which there exists $PD(b,g(t))$ for each $b\in \bf C_s$
with $0<Re (b)$ or $b=0$.
\par $(3)$ Let $S \in Af(X,y)$, $\| S \| _{\mu }<\infty $,
then for each $b \in U:= \{ b': Re \mbox{ }b'>0 \mbox{ or }b'=0 \}$
there is a pseudo-differentiable quasi-invariant measure $\mu $ .}
\par {\bf Proof.} We consider the following additive compact subgroup
$G_T:=\{ x \in X| \| x(j) \| \le p^{k(j)} \mbox{ for each }
j\in {\bf N} \}$ in $X$, where $T=diag \{ d(j)\in K:  \mbox{ }
|d(j)|=p^{-k(j)} \mbox{ for each } j \in {\bf N} \}$ is a compact
diagonal operator.
Then $\mu $ from Theorem 3.15 is quasi-invariant 
relative to the following additive subgroup
$S_T:=G_T+H$, where $H:=sp_{\bf K}\{ e(j): j \in {\bf N} \}$.
The rest of the proof is analogous to that of \S I.4.2.
\par  {\bf 4.3. } Let $X$ be a Banach space over
$\bf K$, $b_0\in {\bf R}$ or $b_0= + \infty $
and suppose that the following conditions are satisfied:
\par $(1)$ $T: X\to X$ is a compact operator with $ker(T)=\{ 0 \}$;
\par $(2)$ a mapping $\tilde F$ from $B({\bf K},0,1)$ to $C_T(X):=\{ U: U 
\in C^1(X,X)$
and $(U'(x)-I)$ is a compact operator for each $x \in X,$
there is $U^{-1}$ satisfying the same conditions as $U \} $ is given;
\par $(3)$ $\tilde F(t)=U_t(x)$ and  $\Phi ^1 U_t(x+h,x)$ are continuous by
$t$, that is, $\tilde F \in C^1(B({\bf K},0,1),C_T(X))$;
\par $(4)$ there is $c>0$ such that
$\| U_t(x)-U_s(x) \| \le \|Tx \|$ for each $x \in X$ and $|t-s| <c$;
\par $(5)$ for each $R>0$ there is a finite-dimensional over $\bf K$ subspace
$H \subset X$ and $c'>0$ such that $\|U_t(x)-U_s(x)\| \le \|Tx\|/R$
for each
$x \in X \ominus H$ and $|t-s| <c'$ with $(3-5)$ satisfying also for
$U^{-1}_t$. 
\par {\bf Theorem.} {\it 
On $X$ there are probability quasi-invariant measures $\mu $ which are
pseudo-differentiable for each $b\in \bf C_s$ with ${\bf R}\ni
Re (b)\le b_0$ relative to a family
$U_t$, where $\mu $ are with values in $\bf K_s$.}
\par {\bf Proof.} From Conditions (2,3) it follows that there is
$c>0$ such that $|det(U_t'(x))|=|det(U_s'(x))|$ in $L(\mu )$
by $x \in X$ and all $|t-s|<c$, where quasi-invariant and pseudo-differentiable
measures $\mu $ on $X$ relative to $S_T$
may be constructed as in the proof of Theorems 3.15 and 4.2.
The final part of the proof is analogous to that of \S I.4.3.
\par {\bf 4.4.} Let $X$ be a Banach space over $\bf K$, $\mu $ be
a probability quasi-invariant measure $\mu : Bco(X) \to \bf K_s$,
that is pseudo-differentiable for a given $b$ with $Re(b)>0$,
$C_b(X)$ be a space of continuous bounded functions
$f: X \to \bf K_s$ with $\| f\| :=
\sup_{x \in X} |f(x)|$. 
\par {\bf Theorem.} {\it For each $a \in J_{\mu }$  and
$f \in C_b(X)$ is defined the following integral:
$$(i) \mbox{ }l(f)=\int _K[\int_X f(x)[\mu (-\lambda a +dx)-
\mu (dx)] g(\lambda ,0,b) v(d \lambda )$$ and there exists
a measure $\nu : Bco(X) \to \bf C_s$
with a bounded variation (for $b \in \bf R$ this $\nu $ is a mapping
from from $Bco(X)$ into $\bf K_s$ such that
$$(ii) \mbox{ } l(f)=\int_X f(x) \nu (dx),$$  where $v$ 
is the Haar measure on $\bf K$ with values in $\bf Q_s$,
moreover, $\nu $ is independent from $f$ and may be dependent on
$a \in J_{\mu }$. We denote $\nu =:\tilde D^b_a \mu $.}
\par {\bf Proof.} From Definition 4.1 and the Lebesgue
theorem it follows that there exists
$\lim_{j \to \infty } \int_{{\bf K}\setminus
B({\bf K},0,p^{-j})}[\int_X (f(x+ \lambda a)-f(x)) g(\lambda ,0,b)
\mu (dx)] v(d \lambda )=l(f)$, that is $(i)$ exists. 
Let $(iii)$ $l_j(
V,f):=$ $\int_{{\bf K}\setminus B({\bf K},0,p^{-j})} 
[\int_V f(x)( \mu (-\lambda a
+dx)- \mu (dx)) g(\lambda ,0,b)] v(d \lambda )$, where $V \in Bco(X)$.
Then due to construction of \S 3.15 for each $c>0$
there exists a compact $V_c \subset X$ with 
$\| X\setminus V_c\| _{\nu _{\lambda }} <c$
for each $|\lambda |>0$, where $\nu _{\lambda }(A):=$
$\int_{{\bf K}\setminus B({\bf K},0,
|\lambda |)} [\mu (-\lambda 'a +A) - \mu (A)] g(\lambda ',0,b)
v(d \lambda ')$ for $A \in Bco(X)$. 
The rest of the proof is analogous to that of \S I.4.4.
\par {\bf 4.5. Theorem.}
{\it Let $X$ be a Banach space over $\bf K$, $|*|=mod_{\bf K}(*)$
with a probability quasi-invariant measure $\mu : Bco(X) \to \bf K_s$ and it
is satisfied Condition $3.16.(i)$, suppose $\mu $ is pseudo-differentiable and
\par $(viii)$ $J_b{\mu } \subset T"J_{\mu }$, $(U_t:$ $ t \in B({\bf K},0,1))$
is a one-parameter family of operators such that
Conditions $3.20(i-vii)$ are satisfied
with the substitution of $J_{\mu }$ onto $J^b_{\mu }$
uniformly by $t \in B({\bf K},0,1)$, $J_{\mu } \supset T'X$, where
$T', T": X \to X$ are compact operators, $ker (T')=ker (T")=0$.
Moreover, suppose that there are sequences 
\par $(ix)$ $[k(i,j)]$ and $[k'(i,j)]$
with $i, j \in \bf N$, $\lim_{i+j \to \infty }k(i,j)=$ $\lim_{i+j\to
\infty }k'(i,j)=- \infty $ and $n \in \bf N$ such that
$|T"_{i,j}-\delta _{i,j}|<|T'_{i,j}-\delta _{i,j}|p^{k(i,j)}$,
$|U_{i,j}-\delta _{i,j}|< |T"_{i,j}-\delta _{i,j}|p^{k'(i,j)}$
and $|(U^{-1})_{i,j}-\delta _{i,j}|<|T"_{i,j}-\delta _{i,j}|p^{k'(i,j)}$
for each $i+j>n$, where $U_{i,j}=\tilde e_iU(e_j)$, $(e_j:$ $j)$ is
orthonormal basis in $X$. Then for each $f \in C_b(X)$ is defined
$$(i)\mbox { }l(f)=\int_{B(K,0,1)}[\int_X f(x)[\mu (U^{-1}_t(dx))-
\mu (dx)]g(t,0,b)v(dt)$$ and there exists a measure $\nu : Bco(X) \to \bf C_s$
with a bounded total variation [particularly, for
$b \in \bf R$ it is such that $\nu : Bco(X)\to \bf K_s$] and
$$(ii)\mbox{ }l(f)=\int_Xf(x)\nu (dx),$$ where $\nu $ is independent
from $f$ and may be dependent on $(U_t:$ $t)$, $\nu =:\tilde D^b_{U_*}\mu $.}
\par {\bf Proof.} From the proof of Theorem 3.20
it follows that there exists a sequence
$U^{(q)}_t$ of polygonal operators converging uniformly by
$t \in B({\bf K},0,1)$ to $U_t$ and equicontinuously
by indices of matrix elements
in $L(\mu )$. Then there exists $\lim_{q\to \infty }
\lim_{j \to \infty }\int_{B({\bf K},0,1)
\setminus B({\bf K},0,p^{-j}}[\int_Xf(U^{-1}_t(x))-f(x)]g(t,0,b)
\mu (dx)]v(dt)$ for each $f \in C_b(X)$. From conditions $(viii, ix),$
the Fubini and Lebesgue theorems it follows that for $\nu _{\lambda }:=
\int_{B({\bf K},0,1)\setminus B({\bf K},0,|\lambda |)}[\mu (U^{-1}_t(A))-
\mu (A)]g(t,0,b)v(dt)$ for $A \in Bco(X)$
for each $c>0$ there exists a compact $V_c \subset X$ and $\delta >0$
such that $\| X \setminus V_c\|<c$.
Indeed, $V_c$ and $\delta >0$ may be chosen due to
pseudo-differentiability of $\mu $, \S \S 2.30, 3.18,
Formula $(i)$, $3.16.(i)$ and due to continuity and boundednessy
(on $B({\bf K},0,1) \ni t$) of $|det$ $U'_t(U^{-1}_t)(x))|_{\bf K}$
satisfying the following conditions $U^{-1}_t(V_c) \subset V_c$
and $\| (X\setminus V_c)
\bigtriangleup (U^{-1}_t(X\setminus V_c))\| _{\mu }=0$ for each
$|t|< \delta $, since $V_c=Y(j)\cap V_c$ are compact for every
$j$. Repeating proofs 3.20 and 4.4 with the use of
Lemma I.2.5 for the family $(U_t:$ $t)$ we get formulas $(i,ii)$.
\section{ Convergence of quasi-invariant and pseudo-differentiable
measures.}
\par {\bf 5.1. Definitions, notes and notations.} Let $S$ be a normal
topological group with the small inductive dimension
$ind (S)=0$, $S'$ be a dense subgroup, suppose their topologies are
$\tau $ and $\tau '$ correspondingly, $\tau ' \supset \tau |_{S'}$.
Let $G$ be an additive Hausdorff left-$R$-module, where $R$ is a
topological ring, ${\sf R } \supset Bco(S)$ be a  
a ring ${\sf R } \supset Bco(S)$ for $\bf K_s$-valued measures,
${\sf M(R,} G)$ be a family of measures with values in $G$,
${\sf L(R,} G,R)$ be a family of quasi-invariant measure 
$\mu : {\sf R} \to G$ with $\rho _{\mu }
(g,x)\times \mu (dx):=\mu ^{g^{-1}}(dx)=:\mu (gdx)$,
$R \times G \to G$ be a continuous left action of $R$ on $G$
such that $\rho _{\mu }(gh,x)=\rho _{\mu}(g,hx)
\rho _{\mu }(h,x)$
for each $g,h \in S'$ and $x \in S$. Particularly,
$1=\rho _{\mu }(g,g^{-1}x)
\rho _{\mu }(g^{-1},x)$, that is, $\rho _{\mu }(g,x) \in R_o$,
where $R_o$ is a multiplicative subgroup of $R$.
Moreover, $zy \in \sf L$ for $z \in R_0$ with $\rho _{z\mu }(g,x)=
z\rho _{\mu }(g,x)z^{-1}$ and $z \ne 0$. We suppose that topological
characters and weights $S$ and $S'$ are countable and each open
$W$ in $S'$ is precompact in $S$.  Let $\bf P"$ be a family of
pseudo-metrics in $G$ generating the initial uniformity such that for each
$c>0$ and $d \in \bf P"$ and $\{U_n\in {\sf R}: n \in {\bf N} \}$ with
$\cap \{ U_n: n \in {\bf N} \}= \{x\}$
there is $m \in \bf N$ such that $d(\mu ^g(U_n), \rho _{\mu }(
g,x)\mu (U_n))<cd(\mu (U_n),0)$ for each $n>m$, in addition,
a limit $\rho $ is independent $\mu $-a.e. on the choice of
$\{U_n: n\}$ for each $x \in S$ and $g \in S'$.
Consider a subring $R' \subset \sf R$,  $R' \supset Bco(S)$
such that $\cup \{A_n: n=1,...,N\} \in R'$
for $A_n \in R'$ with $N \in \bf N$ and $S'R' =R'$. Then
${\sf L(R},G,R;R')
:= \{ (\mu ,\rho _{\mu }(*,*)) \in {\sf L(R},G,R): \mu - \mbox{ }R'-
\mbox{ is regular and for each } s \in S \mbox{ there are } A_n \in R',
n \in {\bf N} \mbox{ with } s=\cap (A_n:n), \{s\} \in R'\}$.
\par  For pseudo-differentiable measures $\mu $ let $S" \subset S'$,
$S"$ be a dense subgroup in $S$, $\tau '|S"$ is not stronger than $\tau "$
on $S"$ and there exists a neighbourhood $\tau " \ni W" \ni e$ in which
are dense elements lying on one-parameter subgroups $(U_t:$
$t\in B({\bf K},0,1)))$. 
We suppose that $\mu $ is induced from the Banach space
$X$ over  $\bf K$ due to a local homeomorphism of neighbourhoods
of $e$ in $S$ and $0$ in $X$ as for the case of groups of
diffeomorphisms \cite{lu1} such that is accomplished Theorem
4.5 for each $U_* \subset S"$ inducing the correspopnding transformations
on $X$. In the following case $S=X$ we consider $S'=J_{\mu }$ and
$S"=J^b_{\mu }$ with $Re(b)>0$ such that $M_{\mu }
\supset J_{\mu } \subset (T_{\mu }X)^{\sim }$, $J^b_{\mu }
\subset (T^{(b)}_{\mu }X)^{\sim }$ with compact operators
$T_{\mu }$ and $T_{\mu }^{(b)}$, $ker (T_{\mu })=ker (T^{(b)}_{\mu })=0$
and norms induced by the Minkowski functional $P_E$ for $E=T_{\mu }
B(X,0,1)$ and $E=T^{(b)}_{\mu }B(X,0,1)$ respectively. 
We suppose furter that for pseudo-differentiable measures
$G$ is equal to $\bf C_s$ $\vee $ $\bf K_s$. We denote
$P({\sf R},G,R,U_*;R'):=[(\mu ,\rho _{\mu }, \eta _{\mu }):$ $(\mu ,
\rho _{\mu }) \in {\sf L(R},G,R;R'), \mu $ is pseudo-differentiable
and $\eta _{\mu }(t,U_*,A) \in L({\bf K},v,{\bf C_s})]$,
where $\eta _{\mu }(t,U_*,A)=j(t)g(t,0,b)[\mu ^h(U^{-1}_t(A)-\mu ^h(A)]$,
$j(t)=1$ for each $t \in \bf K$ for $S=X$; $j(t)=1$ for $t \in B({\bf K},0,1)$,
$j(t)=0$ for $|t|_{\bf K}>1$ for a topological group $S$ 
that is not a Banach space $X$ over $\bf K$,
$v$ is the Haar measure on $\bf K$ with values in $\bf Q_s$,
$(U_t:$ $t \in B({\bf K},0,1))$ is an arbitrary one-parameter subgroup.
On these spaces $\sf L$ (or $P$) the additional conditions are imposed:
\par (a) for each neighbourhood (implying that it is open) $U \ni 0 \in G$
there exists a neighbourhood $S\supset V\ni e$ and a compact subset $V_U$,
$e \in V_U \subset V$, with $\mu (B) \in U$ (or in addition
 $\tilde D^b_{U_*} \mu (B) \in U$) for each $B$, ${\sf R} \ni B
\in Bco(S \setminus V_U)$;
\par  (b) for a given $U$ and a neighbourhood
$R\supset D\ni 0$ there exists a neighbourhood $W$, $S'\supset
W \ni e$, (pseudo)metric $d \in P"$ and $c>0$ such that $\rho _{\mu }
(g,x)-\rho _{\mu }(h,x')
\in D$ (or $\tilde D^b_{U_*}(\mu ^g -\mu ^h)(A) \in U$ for $A \in
Bco(V_U)$ in addition for $P$)
whilst $g, h \in W$, $x,x' \in V_U$, $d(x,x') < c$, where (a,b) is
satisfied for all $(\mu ,\rho _{\mu }) \in \sf L$ (or $(\mu ,\rho _{\mu },
\eta _{\mu }) \in P$) equicontinuously in (a) on $V\ni U_t, U^{-1}_t$ and
in (b) on $W$ and on each
$V_U$ for $\rho _{\mu }(g,x)-\rho _{\mu }(h,x')$
and $\tilde D^b_{U_*}(\mu ^g-\mu ^h)(A)$.
\par These conditions are justified, since 
due to Theorems 3.15, 3.19, 4.3 and 4.5
there exists a subspace $Z"$ dense in $Z'$ such that
for each $\epsilon >0$ and each $\infty >R>0$ there are $r>0$ and $\delta >0$
with $|\rho _{\nu }(g,x)- \rho _{\nu }(h,y) | <\epsilon $ for each
$\| g- h \|_{Z"}+ \| x-y\|_Z<\delta $, $g, h \in B(Z",0,r)$, $x, y\in
B(Z,0,R)$, where $Z"$ is the Banach space over $\bf K$.
For a group of diffeomorphisms of a 
non-Archimedean Banach manifold we have an analogous continuity
of $\rho _{\mu }$ for a subgroup $ G" $ of the entire group $G$
(see \cite{lu1,lu3,lubp2,lu5}).
By $\sf M_o$ we denote a subspace in $\sf M$, satisfying (a).
Henceforth, we imply that $R'$ contains all closed subsets
from $S$ belonging to $\sf R$, where $G$ and $\sf R$ are complete.
\par   For $\mu : Bco(S) \to G$ by $L(S,\mu ,G)$  we denote the completion
of a space of continuous $f: S\to G$ such that $\| f\|_d:=\sup_{
h \in C_b(S,G)}d(\int_Sf(x)h(x)$ $\mu (dx), 0)< \infty $ for each
$d \in P"$, where $C_b(S,G)$ is a space of continuous bounded functions
$h: S\to G$. We suppose that for each sequence $(f_n:$ $n) \subset
L(S,\mu ,G)$ for which  $g \in L(S,\mu ,G)$ exists with
$d(f_n(x),0) \le d(g(x),0)$ for every $d \in P"$, $x$ and $n$,
that $f_n$ converges uniformly on each compact subset $V \subset S$
with $\| V\|_{\mu }>0$ and the following is satisfied: $f \in L(S,\mu ,G)$,
$\lim_n\|f_n-f\|_d=0$ for each $d \in P"$ and $\int_S f(x)\mu (dx)
$ $=\lim_n\int_Sf_n(x)\mu (dx)$. In the case $G=\bf K_s$
it coincides with $L(S,\mu ,{\bf K_s})$, 
hence this supposition is the Lebesgue theorem.
By $Y(v)$ we denote $L({\bf K},v,{\bf C_s})$.
\par   Now we may define topologies and uniformities with the help
of corresponding bases (see below) on ${\sf L}\subset G^{\sf R}
\times R_o^{S'\times S}=:Y$ (or
$P \subset G^{\sf R}\times R_o^{S'\times S}\times
G^{S'\times K\times \sf R}=:Y$, $R_o\subset R\setminus \{0\}$. There are
the natural projections $\pi : \sf L$ $(\vee $ $P)
\to M_o$, $\pi (\mu ,\rho _{\mu }(*,*)$ $(\vee $ ,$\eta _{\mu }))=\mu $,
$\xi: {\sf L}$ $(\vee $ $P) \to R^{S'\times S}$,
$\xi (\mu ,\rho _{\mu },(\vee$ $\eta _{\mu }))=\rho _{\mu }$,
$\zeta : P\to G^{S'\times K\times \sf R}$, $\zeta (\mu ,\rho _{\mu },
\eta _{\mu })=\eta _{\mu }$. Let $\sf H$ be a filter on $\sf L$ or $P$,
$U=U'\times U"$ or $U=U'\times U"\times U"'$, $U'$ and $U"$ be elements
of uniformities on $G$, $R$ and $Y(v)$ correspondingly,
$\tau ' \ni W\ni e$, $\tau \ni V \supset V_{U'} \ni e$, $V_{U'}$ is compact.
By $[\mu ]$ we denote $(\mu ,\rho _{\mu })$ for $\sf L$
or $(\mu ,\rho _{\mu },\eta _{\mu })$ for $P$, $\Omega :=\sf L$ $\vee $
$P$, $[\mu ](A,W,V):=[\mu ^g(A), \rho _{\mu }(g,x),$ $\vee $
$\eta _{\mu ^g}(t,U_*,A)|$ $g \in W, x \in V,$ $\vee $ $t \in K]$.
We consider $\sf A \subset R$, then
$$(1)\mbox{ }{\sf W(A,}W,V_{U'};U):=\{ ([\mu], [\nu ]) \in \Omega ^2|
([\mu ],[\nu ])(A,W,V_{U'}) \subset U \};$$
$$ (2)\mbox{ }{\sf W(S};U):=\{ ([\mu ],[\nu ])\in \Omega ^2|
\{ (B,g,x) : ([\mu ],[\nu ])(B,g,x))\in  U \} \in
{\sf S} \},$$ where  ${\sf S}$ is a filter on ${\sf R} \times
S' \times S^c$, $S^c$ is a family of compact subsets $V' \ni e$.
$$ (3) \mbox{ }{\sf W(F},W,V;U):=\{ ([\mu ],[\nu ]) \in \Omega ^2|
\{ B: ([\mu ],[\nu ])(B,g,x) \in U, g \in W, x \in V \}
\in {\sf F} \},$$ where ${\sf F}$ is a filter on $\sf R$
(compare with $\S $ 2.1 and 4.1\cite{cons});
$$(4) \mbox{ }{\sf W(A,G};U) := \{ ([\mu ],[\nu ]) \in \Omega ^2|
\mbox{ }\{ (g,x): ([\mu ],[\nu])(B,g,x) \in U, B \in {\sf A} \}
\in {\sf G} \},$$ where $\sf G$ is a filter on $S' \times S^c$;
suppose ${\sf U \subset R} \times \tau '_e \times S^c$, $\Phi $
is a family of filters on ${\sf R} \times S' \times S^c$ or
${\sf R}\times S'\times S^c\times Y(v)$ (generated by products of
filters $\Phi_{\sf R} \times \Phi _{S'} \times \Phi _{S^c}$
on the corresponding spaces), ${\sf U'}$ be a uniformity on
$(G,R)$ or $(G,R,Y(v))$,
${\sf F} \subset Y$. A family of finite intersections of sets
${\sf W}(A,U) \cap ({\sf F \times F})$  (see (1)), where
$(A,U) \in {\sf U \times U'}$  (or
$ {\sf W}(F,U) \cap ({\sf F \times F})$ (see (2)), where
$(F,U) \in (\Phi \times {\sf U'}$) generate by the definition
a base of uniformity of $\sf U$-convergence ($\Phi $-convergence
respectively)
on $\sf F$ and generate the corresponding topologies. For these
uniformities are used notations 
$$(i)\mbox{ }{\sf F_U} \mbox{ and }
{\sf F}_ {\Phi };\quad {\sf F}_{{\sf R}
\times W \times V}\mbox{ is for }{\sf F}\mbox{ with the uniformity of uniform
convergence}$$
on ${\sf R} \times W \times V$, where $W \in \tau '_e$,
$V \in S^c$, analogously for the entire space $Y$;
$$(ii)\mbox{ }{\sf F}_A\mbox{ denotes the uniformity (or topology)
of pointwise convergence for }$$
$A \subset {\sf R} \times \tau '_e \times
S^c=:Z$, for $A=Z$ we omit the index (see formula (1)).
Henceforward, we use $\sf H'$ instead of $\mbox{ }\sl H$ in 4.1.24\cite{cons},
that is, ${\sf H'}(A,{\sf \tilde R})$-filter on $\sf R$ generated by the
base $[(L \in {\sf R}:$ $L \subset A \setminus K'):$
$K' \in \sf \tilde R,$
$K' \subset A]$, where ${\sf \tilde R} \subset \sf R$ and $\sf \tilde R$
is closed relative to the finite unions.
\par For example, let $S$ be a locally $\bf K$-convex space, $S'$ be a dense
subspace, $G$ be a locally $\bf L$-convex space, where
$\bf K, L$ are fields, $R=B(G)$ be a space of bounded linear operators on
$G$, $R_o=GL(G)$ be a multiplicative group of invertible linear
operators. Then others possibilities are: $S=X$ be a Banach space 
over $\bf K$, $S'=J_{\mu }$, $S"=J_{\mu }^b$
as above; $S=G(t)$, $S'\supset S"$ are dense subgroups,
$G=R$ be the field $\bf K_s$ ( $s \ne p$), $M$ be an
analytic Banach manifold over $\bf K \supset
Q_p$ (see \cite{lu1}). The rest of the necessary standard definitions
are recalled further when they are used.
\par {\bf 5.2. Lemma.} {\it Let $\sf R$ be a quasi-$\delta $-ring
with the weakest uniformity in which each $\mu  \in \sf M$ is uniformly
continuous and $\Phi \subset \hat \Phi _C({\sf R}, S' \times S^c)$. Then
${\sf L(R}, G,R,R')_{\Phi }$ (or $P({\sf R},G,R,U_*;R')_{\Phi }$) )
is a topological space on which $R_o$ acts continuously from the right.}
\par {\bf Proof.} It is analogous to that of \S I.5.2
using
Definition 4.1 for pseudo-differentiable $f$.
\par {\bf 5.3. Proposition.} {\it (1). Let $\sf T$ be a
$\hat \Phi _4$-filter on ${\sf M_o(R,}G;R')$, $\{ A_n\}$ be a disjoint
$\Theta ({\sf R})$-sequence, $\Sigma $ be the elementary filter on $\sf R$
generated by $\{A_n:n \in {\bf N} \}$ and $\phi : {\sf M_o \times R}
\to G$ with $\phi (\mu ,A)=\mu (A)$. Then $\phi ({\sf T} \times \Sigma )$
converges to $0$.  (2). Moreover, let $\sf U$ be a base of neighbourhoods
of $e \in S'$, $\phi: {\sf L} \to G\times R$, $\phi (\mu ,A,g):=(\mu ^g(A),
\rho _{\mu }(g,x))$, where $x \in A$. Then $(0,1) \in
\lim \phi ({\sf T} \times \Sigma \times {\sf U})$.  (3). If $T$ is a
$\Phi _4$-filter on $P({\sf R},G,R,U_*;R')$, $\psi (\mu ,B,g,t,U_*)=
[\mu (B); \rho _{\mu }(g,x); \eta _{\mu ^g}(t,U_*,B)]$, then
$(0,1,0) \in \lim \psi ({\sf T}\times \Sigma \times {\sf U})$ for each
given $U_* \in S"$, where $\Sigma $ and $\sf U$ as in (1,2).}
\par {\bf Proof.} The proof is analogous to that of
\S I.5.3 with the use of the Lebesgue
convergence theorem.
\par {\bf 5.4. Proposition.} {\it Let $\sf H$ be a
$\hat \Phi _4 $-filter
on $\sf L$ (or $P$) with the topology $\sf F$ (see 5.1.(ii)),
$A \in \sf R$,
$\tau _G \ni U \ni 0$, ${\sf H'}(A,R') \in \Psi _f({\sf R})$. Then
there are $L \in \sf H$, $\tilde K \in R'$ and an element of the
uniformity $\sf U $ for ${\sf L}_{R'}$ or $P_{R'}$ such that
$\tilde {\bf K} \subset A$,
$L=[(\mu ,\rho _{\mu }(g,x)):$
$M:=\pi _{\sf M_o}(L) \ni \mu , \pi _{\tau '_e}(L)=:W \ni g $ (or
$(\mu ,\rho _{\mu },\eta _{\mu }(*,*,U_*))$ and additionally
$\tilde D^b_{U_*}\mu
=PD(b,\eta _{\mu }))]$, $e \in W \in \tau '$,
$\mu ^g(B)-\nu ^{h}(C) \in U$ (or in addition $(\tilde D^b_{U_*}
\mu ^g(B)) -(\tilde D^b_{U_*}\nu ^h(C)) \in U$)
for $\tilde K \subset B \subset A$, $\tilde K \subset
C \subset A$ for each $([\mu ],[\nu ]) \in \bar {\bar L}^2
\cap {\sf U}$, where $\bar {\bar L}:=cl(L,{\sf L}_{R'})$ (or
$cl(L,P_{R'})$), $\pi _{\sf M_o}$ is a projector from
$L$ into $\sf M_o$.}
\par {\bf Proof.} Repeating the proof of \S I.5.4 we get
$\mu ^g(B)-\mu (B) \in U'$, $\nu ^h(C)-
\nu (C) \in U'$ and for $3U' \subset U$ we get 5.4 for $\sf L$.
From Theorems 4.4 and 4.5, \S 5.1, the Egorov conditions and the Lebesgue
theorem we get 5.4 for $P$, since $\mu $ are probability measures
and ${\sf L}_{R'}$ (or $P_{R'}$) correspond to uniformity
from \S 5.1.(ii) with $A=R'\times
\tau _e'\times S^c$. Indeed, $\mu ^g(A)-\nu ^h(A)=
(\mu ^g(A)-\mu ^g(V_{U'}))+(\mu ^g(V_{U'})-\nu ^h(V_{U'}))+(\nu ^h(V_{U'})-
\nu ^h(A))$, $\mu ^g(A)=\int_A\rho _{\mu }(g,x)\mu (dx)$ for
each $A \in Bco(S)$, for each $\tau _G\ni U'\ni 0$ there exists a
compact subset $V_U' \subset A$
with $\mu ^g(B)
\in U'$ for each $B \in Bf(A\setminus V_{U'})\cap Bco(S)$ and the same
for $\nu ^h$ (due to the condition in \S 5.1 that $R'$ contains $Bco(S)$).
At first we can consider $A \in Bco(S)$, then use
$R'$-regularity of measures and $\sigma R'\supset Bco(S)$.
From the separability of $S$, $S'$ and the equality of their
topological weights to $\aleph _0$, restrictions 5.1.(a,b) it follows
that there exists a sequence of partitions $Z_n=[(x_m,A_m):$
$m, x_m \in A_m]$
for each $A \in Bco(S)$, $A_i\cap A_j=\emptyset $
for each $i \ne j$, $\bigcup_mA_m=A$, $A_m \in Bco(S)$,
such that $\lim_{n\to \infty }(\mu ^g(A)-\sum_j\rho _{\mu }(g,x_j)\mu
(A_j))=0$ and the same for $\nu $, moreover, for $V_{U'}$ each
$Z_n$ may be chosen finite. Then there exists $W \in \tau '_e$ with
$W\times
(S\setminus V^2)\subset (S\setminus V)$, $\tau _e\ni V \subset V^2$,
$\nu ^g(B)$ and
$\mu ^g(B) \in U'$ for each
$B \in Bf(S\setminus V^2)\cap Bco(S)$ (for $G=
\bf K_s$ respectively) and $g \in W$ (see 5.1.(a)). Then from $A=[A\cap
(S\setminus V^2)]\cup [A\cap V^2]$ and the existence of compact
$V'_{U'}
\subset V$ with $\mu (E) \in U'$ for each $E \in Bf(V \setminus V'_{U'})
\cap Bco(S)$ and the same for $\nu $, moreover, $(V'_{U'})^2$ is also
compact,
it follows that $\mu ^g(B)-\nu ^h(C)\in U$ for $9U' \subset U$, since
$R' \supset Bco(S)$, where $W$ satisfies the following condition
$\mu ^g(V'_{U'})-\nu ^h(V'_{U'}) \in U'$ for $V'_{U'} \subset V^2$
due to \S 5.1.(b), $\mu (B)-\nu (C) \in U'$, $WV'_{U'} \subset (V'_{U'})^2$
due to precompactness of $W$ in $S$. Since pseudo-differentiable
measures are also quasi-invariant, hence for them 5.4 is true.
\par  Now let $[\mu ] \in \lim H$, $A \in Bco(S)$, then $\eta _{\mu }
\in \lim \zeta (H)$ in $Y(v)$ and there exists a sequence
$\eta _{\mu _n}$ such that $\int_{\bf K}\eta _{\mu _n}(\lambda ,U_*,A)
v(d \lambda )=\tilde D^b_{U_*}\mu _n(A)$ due to \S \S 4.4 or 5.1 and
$\lim _{n \to \infty }\tilde D^b_{U_*}\mu _n(A)=$ $\int_{\bf K}\eta _{\mu }
(\lambda ,U_*,A)v(d \lambda )=:\kappa (A)$ due to the Lebesgue theorem.
From $\eta _{\mu }(\lambda ,U_*,A\cup B)=\eta (\lambda ,U_*,A)+
\eta (\lambda ,U_*,B)$ for $A\cap B=\emptyset $, $B \in Bco(S)$
it follows that $\nu (A)$ is the
measure on $Bco(S)$, moreover,
$\kappa (A)=\tilde D^b_{U_*}\mu (A)$. Since $\mu ^g(A)=\int_A\rho _{\mu }
(g,x)\mu (dx)$ for $A \in Bco(S)$ for $g\in S'$, then $\eta _{\mu ^g}
(\lambda ,U_*,A)=j(\lambda )g(\lambda ,0,b)[\mu ^g(A)-\mu ^g(U^{-1}
_{\lambda }A)]=$ $j(\lambda )g(\lambda ,0,b)\int_A\rho _{\mu }
(g,x)[\mu (dx)-\mu ^{U_{\lambda }}(dx)]$ and in view of the Fubini theorem
there exists $\tilde D^b_{U_*}\mu ^g(A)=\int_A[\int_K\rho _{\mu }
(g,x)$ $j(\lambda )g(\lambda ,0,b)[\mu (dx)-\mu ^{U_{\lambda }}(dx)]v(d
\lambda )$, where $j(t)=1$ for $S=X$ and $j(t)$ is the characteristic
function of $B({\bf K},0,1)$ for $S$ that is not the Banach space $X$. 
Then $\mu $-a.e.
$\tilde D^b_{U_*}\mu ^g(dx)/ \tilde D^b_{U_*}\mu (dx)$ coincides with
$\rho _{\mu }(g,x)$ due to 5.1.(a,b), hence, $(\tilde D^b_{U_*}
\mu ^g, \rho _{\mu ^g})$ generate the $\Phi _4$-filter in $L$ arising
from the $\hat \Phi _4$-filter in $P$. Then we estimate
$\tilde D^b_{U_*}(\mu ^g-\nu ^h)(A)$ as above $\mu ^g(A)-\nu ^h(A)$.
Therefore, we find for the $\Phi _4$-filter corresponding $L$,
since there exists $\delta >0$
such that $U_{\lambda }\in W$ for each $|\lambda |<\delta $. For
$\Phi _4$-filter we use the corresponding finite intersections
$W_1\cap ...\cap W_n=W$, where $W_j$ correspond to the
$\Phi _4$-filters $H_j$.
\par {\bf Note.} The formulations and proofs of \S \S 5.5-5.10 
(see Part I) are quite
analogous for real-valued and $\bf K_s$-valued measures
due to preceding results.

\par Address: Theoretical Department,
\par Institute of General Physics,
\par Russian Academy of Sciences,
\par Str. Vavilov 38, Moscow, 119991 GSP-1, Russia
\end{document}